\documentclass[10pt, a4paper, oneside]{article}
\usepackage[italian,english]{babel}
\usepackage[a4paper,top=30mm,bottom=40mm,inner=30mm,outer=20mm]{geometry}
\usepackage{listings}
\usepackage{color}
\usepackage{setspace}
\usepackage{hyperref}
\usepackage{acro}
\usepackage{amsmath}
\usepackage{amsthm}
\usepackage{amssymb}
\usepackage{mathtools}
\usepackage{graphicx}
\usepackage{geometry}
\usepackage[font=footnotesize,labelfont=bf]{caption}
\usepackage[font=footnotesize]{subcaption}
\usepackage{bbm} 
\usepackage{bm} 
\usepackage{mdframed} 
\usepackage{comment}
\usepackage{siunitx}
\usepackage{hyperref}
\usepackage{xurl} 
\usepackage{tikz}
\usetikzlibrary{arrows}
\usetikzlibrary{decorations.pathreplacing}
\usepackage{todonotes}
\usepackage{enumitem}
\usepackage{array}
\usepackage{xcolor}

\usepackage{cancel}
\usepackage{tikz}
\usepackage{color}
\usetikzlibrary{calc,trees,positioning,arrows,chains,shapes.geometric,%
    decorations.pathreplacing,decorations.pathmorphing,shapes,%
    matrix,shapes.symbols}
\include{acronyms}
\hypersetup{
    colorlinks,
    citecolor=black,
    filecolor=black,
    linkcolor=black,
    urlcolor=black
}

\definecolor{DarkerGreen}{RGB}{0,179,45}

\include{python_highlighting}

\newtheorem{exmp}{Example}[section]

\tikzset{
>=stealth',
  punktchain/.style={
    rectangle, 
    rounded corners, 
    draw=black, very thick,
    text width=10em, 
    minimum height=3em, 
    text centered, 
    on chain},
  line/.style={draw, thick, <-},
  element/.style={
    tape,
    top color=white,
    bottom color=blue!50!black!60!,
    minimum width=8em,
    draw=blue!40!black!90, very thick,
    text width=10em, 
    minimum height=3.5em, 
    text centered, 
    on chain},
  every join/.style={->, thick,shorten >=1pt},
  decoration={brace},
  tuborg/.style={decorate},
  tubnode/.style={midway, right=2pt},
}

\allowdisplaybreaks[3]

\usepackage[toc,page]{appendix} 

\newcounter{for}[section]

\newtheorem{itlemma}{Lemma}[section]
\newtheorem{itproposition}[itlemma]{Proposition}
\newtheorem{itfact}[itlemma]{Fact}
\newtheorem{theorem}[itlemma]{Theorem}
\newtheorem{itcorollary}[itlemma]{Corollary}
\newtheorem{itremark}[itlemma]{Remark}
\newtheorem{itremarks}[itlemma]{Remarks}
\newtheorem{itdefinition}[itlemma]{Definition}
\newtheorem{itexample}[itlemma]{Example}

\newenvironment{fact}{\begin{itfact}\rm}{\end{itfact}}
\newenvironment{claim}{\begin{itclaim}\rm}{\end{itclaim}}
\newenvironment{lemma}{\begin{itlemma}}{\end{itlemma}}
\newenvironment{remark}{\begin{itremark}\rm}{\end{itremark}}
\newenvironment{remarks}{\begin{itremarks} \rm}{\end{itremarks}}
\newenvironment{corollary}{\begin{itcorollary}}{\end{itcorollary}}
\newenvironment{proposition}{\begin{itproposition}}{\end{itproposition}}
\newenvironment{definition}{\begin{itdefinition}\rm}{\end{itdefinition}}
\newenvironment{example}{\begin{itexample}\rm}{\end{itexample}}
\newcommand{\be}[1]{\addtocounter{for}{1} \begin{equation}\label{#1}}
\newcommand{\ee}{\end{equation}}
\newcommand{\bl}[1]{\begin{lemma}\label{#1}}
\newcommand{\br}[1]{\begin{remark}\label{#1}}
\newcommand{\brs}[1]{\begin{remarks}\label{#1}}
\newcommand{\bt}[1]{\begin{theorem}\label{#1}}
\newcommand{\bd}[1]{\begin{definition}\label{#1}}
\newcommand{\bp}[1]{\begin{proposition}\label{#1}}
\newcommand{\bfact}[1]{\begin{fact}\label{#1}}
\newcommand{\bc}[1]{\begin{corollary}\label{#1}}
\newcommand{\bex}[1]{\begin{example}\label{#1}}
\newcommand{\ec}{\end{corollary}}
\newcommand{\efact}{\end{fact}}
\newcommand{\eex}{\end{example}}
\newcommand{\el}{\end{lemma}}
\newcommand{\er}{\end{remark}}
\newcommand{\ers}{\end{remarks}}
\newcommand{\et}{\end{theorem}}
\newcommand{\ed}{\end{definition}}
\newcommand{\ep}{\end{proposition}}
\newcommand{\epr}{\end{proof}}
\newcommand{\bpr}{\begin{proof}}
\newcommand{\bcl}[1]{\begin{claim}\label{#1}}
\newcommand{\ecl}{\end{claim}}

\newcommand{\ind}{{\bf{1}}}

\newcommand{\ecs}{\end{corollary}}
\newcommand{\eers}{\end{exercise}}
\newcommand{\eexs}{\end{example}}
\newcommand{\eems}{\end{example}}
\newcommand{\els}{\end{lemma}}
\newcommand{\eles}{\end{lemmaex}}
\newcommand{\ets}{\end{theorem}}
\newcommand{\eds}{\end{definition}}
\newcommand{\eps}{\end{proposition}}

\newcommand{\bi}{\begin{itemize}}
\newcommand{\ei}{\end{itemize}}
\newcommand{\ben}{\begin{enumerate}}
\newcommand{\een}{\end{enumerate}}

\def\vbar{\mathchoice{\vrule height6.3ptdepth-.5ptwidth.8pt\kern-.8pt}
   {\vrule height6.3ptdepth-.5ptwidth.8pt\kern-.8pt}
   {\vrule height4.1ptdepth-.35ptwidth.6pt\kern-.6pt}
   {\vrule height3.1ptdepth-.25ptwidth.5pt\kern-.5pt}}
\def\fudge{\mathchoice{}{}{\mkern.5mu}{\mkern.8mu}}
\def\bbc#1#2{{\rm \mkern#2mu\vbar\mkern-#2mu#1}}
\def\bbb#1{{\rm I\mkern-3.5mu #1}}
\def\bba#1#2{{\rm #1\mkern-#2mu\fudge #1}}
\def\bb#1{{\count4=`#1 \advance\count4by-64 \ifcase\count4\or\bba A{11.5}\or
   \bbb B\or\bbc C{5}\or\bbb D\or\bbb E\or\bbb F \or\bbc G{5}\or\bbb H\or
   \bbb I\or\bbc J{3}\or\bbb K\or\bbb L \or\bbb M\or\bbb N\or\bbc O{5} \or
   \bbb P\or\bbc Q{5}\or\bbb R\or\bbc S{4.2}\or\bba T{10.5}\or\bbc U{5}\or
   \bba V{12}\or\bba W{16.5}\or\bba X{11}\or\bba Y{11.7}\or\bba Z{7.5}\fi}}

\def \R {{\mathbb R}}

\def \ra {\rightarrow }

\def \o {\omega}
\def \a{\alpha}

\def \s{\sigma}

\def \P{{\mathbb{P}}}

\def\e{\varepsilon}

\def\L{\Lambda}
\def\l{\lambda}
\def\E{{\mathbb{E}}}

\def\1{{\bf 1}}
\def\e{\varepsilon}

\def\ec{\`e }

\begin{document}

\title{Noise-induced oscillations  \\ for 
 the mean-field Dissipative Contact Process}
\author{Paolo Dai~Pra, Elisa Marini}

\maketitle

\begin{abstract}
We study a dissipative version of the contact process, with mean-field interaction, which admits a simple epidemiological interpretation. The propagation of chaos and the corresponding normal fluctuations reveal that the noise present in the finite-size system induces oscillations with a nearly deterministic period and a randomly varying amplitude. This is reminiscent of the emergence of pandemic waves in real epidemics. 

\end{abstract}

\section{Introduction}

The contact process is the microscopic counterpart of one of the most basic epidemiological models, the SIS model. Suppose the individuals of a population are placed at the vertices of a graph $(V,E)$; the individual at the vertex $i \in V$ can be either {\em susceptible} ($x_i = 0$) or {\em non susceptible} ($x_i = 1$). The configuration $x = (x_i)_{i \in V}$ evolves as a continuous time Markov chain with the following rates:
\bi
\item
each non susceptible individual becomes susceptible with rate $1$;
\item
each susceptible individual becomes non susceptible with a rate equal to a given constant $\l$ times the fraction of non susceptible neighbors.
\ei 
This definition suffices in finite graphs, while further care is needed in countable graphs (see \cite{L13}). In finite graphs the chain is absorbed in the null state ($x_i \equiv 0$) in finite time (we say that the epidemics {\em dies out}). The {\em time to absorption} could however vary from a {\em fast absorption} (at most polynomial in $|V|$) to a slow absorption (exponential in $|V|$). For many (sequences of) graphs the transition from fast to slow absorption occurs as the {\em infection constant} $\l$ crosses a critical value $\overline{\l}$. The simplest  nontrivial example is the case in which $(V,E)$ is the {\em complete graph} of $N$ vertices:  $V = \{1,2,\ldots,N\}$, $E = V \times V$. The study of this model reduces to the analysis of the one-dimensional Markov process
\[
m_N := \frac{1}{N} \sum_{i=1}^N x_i,
\]
and the critical infection constant is $\overline{\l} = 1$. Infinite countable graphs require a more detailed analysis, leading to results in several directions, including ergodicity, convergence, local and global survival of the epidemics; such rich behavior could imply the existence of more than one critical value for the infection constant.

In this paper we introduce a modification of the contact process. The state $x_i$ of each individual takes values in the interval $[0,1]$, and is interpreted as her {\em viral load}; we say that an individual is susceptible if $x_i=0$. The dynamics goes as follows:
\bi
\item
each non susceptible individual becomes susceptible with rate $r>0$;
\item
each susceptible individual becomes non susceptible with a rate equal to $\l>0$ times the arithmetic mean of the viral loads of her neighbors, and her viral load jumps to the value $1$; 
\item
between jumps, the viral load decays exponentially with rate $\a>0$.

\ei 
 Note that the model is overparametrized, as we could rescale the time to have $r=1$. We keep however all parameters as they have a useful interpretation: $r^{-1}$ is the time scale at which infected (i.e. non susceptible) individuals lose their immunity, while $\a^{-1}$ is the time scale at which infected individuals remain contagious. In many real epidemics $\a^{-1} \ll r^{-1}$: for a relatively long time non susceptible individuals are {\em immune}, and do not contribute to the propagation of the disease.  Note that many of the useful mathematical properties of the contact process are lost; for instance, monotonicity breaks due to the presence of immune individuals that block contagion.
 
In this paper we study this modified contact process, that we call {\em dissipative}, on the complete graph. Unlike the corresponding classical contact process, the dissipative version does not admit a finite dimensional reduction: any finite-dimensional functional of the empirical measure is non Markovian for large $N$. We particularly deal with the thermodynamic limit of the process ($N \ra +\infty$), and we obtain a law of large numbers (propagation of chaos) and a central limit theorem. Despite the lack of finite dimensional reduction, both the {\em limit process} corresponding to the law of large numbers and the limit {\em fluctuation process} corresponding to the central limit theorem are two-dimensional. Under suitable conditions on the parameters, including the case $\a^{-1} \ll r^{-1}$, the limit process has a unique stable fixed point, which is approached by damped oscillations. The limit fluctuation process reveals, via a Fourier analysis, that noise induces persistent oscillations, despite of the damping exhibited by the limit process. The nature of these oscillations is then investigated by letting the parameters $\a$ and $\l$ to diverge (slowly) with $N$: if properly rescaled, the dynamics of the fraction of non susceptible individuals and of the mean viral load converges to a harmonic oscillator affected by noise.
In their epidemiological interpretation these results suggest that pandemic waves are induced by finite-size effects, allowing the random noise to excite frequencies close to a characteristic value. The resulting motion is however not strictly periodic: though the period is nearly constant, the amplitude of each oscillation is significantly affected by noise, in agreement with what observed in real epidemics.

This paper is organized as follows. Section \ref{sec:micro_dyn_and_prop_ch} deals with the derivation of the dynamics of a representative individual in the limit as $N \ra +\infty$. Detailed properties of this limit dynamics are studied in Section \ref{sec:macro_dyn}. In Section \ref{sec:normal_fluct} we study the dynamics of normal fluctuations around the limit dynamics. In Section \ref{sec:rescale} we scale with $N$ the parameters $\a$ and $\l$ of the model, and study the corresponding limit dynamics. Unless otherwise stated, all proof are collected in Section \ref{sec:proofs}.


\section{Microscopic dynamics and propagation of chaos}\label{sec:micro_dyn_and_prop_ch}

\begin{figure}[!hb]
\centering
\tikzset{every picture/.style={line width=0.75pt}} 

\begin{tikzpicture}[x=0.75pt,y=0.75pt,yscale=-0.75,xscale=0.75]
\draw   (59,150) .. controls (59,127.36) and (77.36,109) .. (100,109) .. controls (122.64,109) and (141,127.36) .. (141,150) .. controls (141,172.64) and (122.64,191) .. (100,191) .. controls (77.36,191) and (59,172.64) .. (59,150) -- cycle ;
\draw    (101,109) .. controls (140.2,79.6) and (208.35,79.56) .. (248.57,107.27) ;
\draw [shift={(251,109)}, rotate = 216.44] [fill={rgb, 255:red, 0; green, 0; blue, 0 }  ][line width=0.08]  [draw opacity=0] (8.93,-4.29) -- (0,0) -- (8.93,4.29) -- cycle    ;
\draw    (103.91,193.24) .. controls (146.08,221.55) and (357.24,220.81) .. (400,191) ;
\draw [shift={(101,191)}, rotate = 41.84] [fill={rgb, 255:red, 0; green, 0; blue, 0 }  ][line width=0.08]  [draw opacity=0] (8.93,-4.29) -- (0,0) -- (8.93,4.29) -- cycle    ;
\draw   (210,150) .. controls (210,127.36) and (228.36,109) .. (251,109) .. controls (273.64,109) and (292,127.36) .. (292,150) .. controls (292,172.64) and (273.64,191) .. (251,191) .. controls (228.36,191) and (210,172.64) .. (210,150) -- cycle ;
\draw   (359,150) .. controls (359,127.36) and (377.36,109) .. (400,109) .. controls (422.64,109) and (441,127.36) .. (441,150) .. controls (441,172.64) and (422.64,191) .. (400,191) .. controls (377.36,191) and (359,172.64) .. (359,150) -- cycle ;
\draw    (292,150) .. controls (331.2,120.6) and (319.5,177.64) .. (356.66,151.69) ;
\draw [shift={(359,150)}, rotate = 143.13] [fill={rgb, 255:red, 0; green, 0; blue, 0 }  ][line width=0.08]  [draw opacity=0] (8.93,-4.29) -- (0,0) -- (8.93,4.29) -- cycle    ;
\draw [shift={(328.1,152.81)}, rotate = 227.43] [fill={rgb, 255:red, 0; green, 0; blue, 0 }  ][line width=0.08]  [draw opacity=0] (8.93,-4.29) -- (0,0) -- (8.93,4.29) -- cycle    ;

\draw (93,142) node [anchor=north west][inner sep=0.75pt]  [font=\small] [align=center] {$\displaystyle  0$};
\draw (243,142) node [anchor=north west][inner sep=0.75pt]  [font=\small] [align=left] {$\displaystyle  1$};
\draw (385,139) node [anchor=north west][inner sep=0.75pt]  [font=\small] [align=left] {$\displaystyle  e^{-\alpha t}$};
\draw (152,30) node [anchor=north west][inner sep=0.75pt]  [font=\small] [align=left] {$\displaystyle \frac{\l}{N}\sum _{j =1}^N x_j$};
\draw (245,220) node [anchor=north west][inner sep=0.75pt]  [font=\small] [align=left] {$\displaystyle r$};
\draw (324,126) node [anchor=north west][inner sep=0.75pt]  [font=\small] [align=left] {$\displaystyle \alpha $};

\end{tikzpicture}
\caption{\footnotesize{An individual $i$ who is initially susceptible ($x_i=0$), can become infectious with rate $\frac{\l}{N}\sum_{j=1}^{N}x_j$ (mean-field interaction). When this happens, her viral load $x_i$  jumps to the value $1$. When infectious, individual $i$ can recover with rate $r$ and, in this case, she becomes immediately susceptible again. In the meantime, the viral load of the individual decreases deterministically in time with rate $\alpha$. Notice that this dynamics is invariant under permutations of the particles. This model is a piecewise deterministic Markov process, as the $x_i$s have a deterministic continuous dynamics with discontinuities determined by random jump events. 
}}
\label{fig:tikz_modello}
\end{figure}
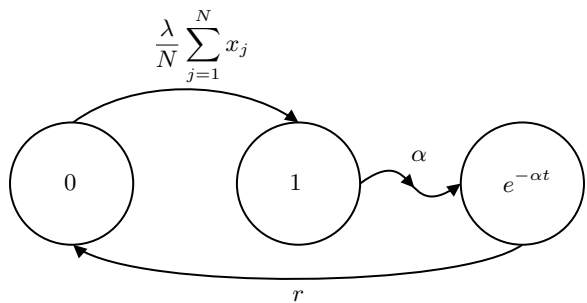

We consider a population of $N$ individuals whose states $(x_i)_{i=1}^N$ take values in $[0,1]$.  
The microscopic dynamics of the system is sketched in Figure \ref{fig:tikz_modello}.

We define $m_N(t)$ and $v_N(t)$ to be the average number of infected individuals and the average viral load of the $N$-particle system at time $t$ respectively, namely, 
\begin{equation}
    m_N(t) \coloneqq \frac{1}{N}\sum_{i=1}^N \ind_{\{x_i(t) >0 \}} = \int_0^1 \ind_{\{y>0\}} \rho^N_t(dy)
 \end{equation}
and 
\begin{equation}
    v_N(t) \coloneqq \frac{1}{N}\sum_{i=1}^N x_i(t) = \int_0^1 y \rho^N_t(dy)
\end{equation}
where $\rho^N_t$ is the empirical measure of $\left(x_i\right)_{i=1}^N$ at time $t$, i.e., the measure-valued process
\begin{equation}
	\rho^N_t \coloneqq \frac{1}{N}\sum_{i=1}^N \delta_{x_i(t)}
\end{equation}
The state $\mathcal X^N \coloneqq \left(x_i\right)_{i=1}^N$, taking values on $[0,1]^N$, is a Markov process whose infinitesimal generator is the closure of the operator $\mathcal{L}^N$ acting on continuously differentiable functions $\phi \in C^1\left([0,1]^N\right)$ as 
\begin{equation}\label{eq:micro:dyn:generator}
    \begin{split}
        \mathcal{L}^{N}\phi\left(\mathcal X^N\right) & = \sum_{k=1}^N -\alpha x_k \partial_{x_k} \phi\left(\mathcal X^N\right) \\
        & + \sum_{k=1}^N  \ind_{\{x_k(t) =0 \}} \l v_N \left[\phi\left(\mathcal X^{N,\uparrow,k} \right) - \phi\left(\mathcal X^N\right)\right]\\
        & + \sum_{k=1}^N   \ind_{\{x_k(t) >0 \}} r \left[ \phi\left( \mathcal  X^{N,\downarrow,k} \right) -\phi\left(\mathcal  X^N\right)\right]
    \end{split}
\end{equation}
where $\mathcal X^{N,\uparrow,k}$ denotes a state equal to $\mathcal X^N$, except for the $k$-th coordinate, which is $1$, and $\mathcal X^{N,\downarrow,k}$ denotes a state equal to $\mathcal X^N$, except for the $k$-th coordinate, which is $0$.  

Equivalently, the state $x_i(t)$ of each individual evolves according to 
\begin{equation}\label{eq:micro:dyn}
\begin{split}
     x_i(t) & = x_i(0) - \int_0^t \alpha x_i(s) ds +  \int_{[0,t]\times [0,+\infty)} \pmb{1}_{[0, \ind_{\{x_i(s^-) =0 \}}  \l v_N(s^-)]}(u) N^\uparrow_i(ds,du) \\
     &  - \int_{[0,t]\times [0,+\infty)} x_i(s^-) \pmb{1}_{[0, \ind_{\{x_i(s^-) >0 \}}(s^-) r]}(u) N^\downarrow_i(ds,du)
\end{split}
\end{equation}
Here $(N^\uparrow_i )_{i=1}^N$ and $(N^\downarrow_i )_{i=1}^N$ are families of independent Poisson random measures on $[0,+\infty) \times [0,+\infty)$ with intensity measure the Lebesgue measure $dt \times du$. Moreover, $ N^\uparrow_i$ and $ N^\downarrow_j$ are independent for all $i,j \in \left\{1,\dots, N\right\}$. 
Notice that integrating over $[0,t]$ makes the processes $x_i$ càdlàg. Also, here $f(t^-) \coloneqq \lim_{s\to t^-} f(s)$ is standard notation for the left limit of a function $f$ at $t$. Strong existence and uniqueness for system \eqref{eq:micro:dyn} can be proven as in \cite{G92}. 

\vspace{20pt}

We will prove that the evolution of a representative component in the limit as $N \ra +\infty$ is given by the stochastic differential equation
\begin{equation}\label{eq:prop:chaos}
    \begin{split}
         \bar x(t) & = \bar x(0) - \int_0^t \alpha \bar x(s) ds +  \int_{[0,t] \times [0,+\infty)} \pmb{1}_{[0, \ind_{\{\bar x(s^-)=0\}}\l v(s^-)]}(u) N^\uparrow(ds,du) \\
         &  - \int_{[0,t] \times [0,+\infty)} \bar x(s^-) \pmb{1}_{[0, \ind_{\{\bar x(s^-)>0\}} r ]}(u) N^\downarrow(ds,du)
    \end{split}
\end{equation}
where $N^\uparrow$ and $N^\downarrow$ are independent Poisson random measures on $[0,+\infty)\times [0,+\infty)$, both having intensity measure $dt\times du$, and $v(t) \coloneqq  \int_{[0,+\infty)} u \rho_{\bar{x}(t)}(du)$ is the first moment of the law $\rho_{\bar{x}}$ of the variable $\bar x$ at time $t$.  

Notice that Eq. \eqref{eq:prop:chaos} is a nonlinear, also called McKean-Vlasov, SDE, as the law of its solution appears as an argument of its coefficients through its first moment $v$.  
Existence and pathwise uniqueness of a strong solution to system \eqref{eq:prop:chaos} follow straightforwardly from \cite{DPMarsiglia} and \cite{G92}.

\vspace{20pt}

\begin{theorem}[\textbf{Propagation of chaos}]\label{thm:prop:chaos}
Suppose $\left(x_i(0)\right)_{i = 1}^N$ are independent identically distributed random variables with values in $[0,1]$ and law $\mu_0$.   
Denote by  $(\left(x_i(t)\right)_{t\geq 0})_{i=1}^N$ the corresponding solution to system \eqref{eq:micro:dyn}. 
Also, consider $N$ independent copies of the solution to Eq. \eqref{eq:prop:chaos}, $(\left(\bar x_i(t)\right)_{t\geq 0})_{i=1}^N$, with the same Poisson random measures $N^\uparrow_i$, $N^\downarrow_i$, $i=1,\dots,N$ and the same initial conditions of system \eqref{eq:micro:dyn}.  

Then, for all $i=1,\dots,N$ and for all $T>0$, 
\begin{equation}\label{eq:bound:prop:chaos}
    \E\left[\sup_{t\in[0,T]} \left\{\left|x_i(t)-\bar{x}_i(t) \right|+\left| \ind_{\{x_i(t) >0 \}}-\ind_{\{\bar x_i(t)>0\}} \right|\right\}\right] \leq \frac{A_T}{\sqrt{N}}
\end{equation}
where 
\begin{equation}\label{C_T}
A_T = \l T e^{(\a+3\l + 2r)T}.
\end{equation}
\end{theorem} 

The proof of Theorem \ref{thm:prop:chaos} is a standard application of coupling arguments and Gromwall Lemma. We however include it in this paper since the specific form of the constant $A_T$ will be useful in Section \ref{sec:rescale}. 

Theorem \ref{thm:prop:chaos}, together with the fact that the i.i.d. processes $\bar{x}_i$ satisfy a standard Law of Large Numbers, implies the following result.  
\begin{corollary}\label{cor:LLN}
For every $T>0$
\begin{equation}\label{eq:correct:rescaling:for:the:fluctuations:uniform}
    \begin{split}
        \E\left[\sup_{t\in[0,T]} \left| m_N(t) - \P\left[\bar{x}(t) >0\right]\right|\right] & \leq \frac{B_T}{\sqrt{N}} \qquad \qquad \E\left[\sup_{t \in [0,T]} \left| v_N(t) - \E\left[\bar{x}(t)\right] \right| \right] \leq \frac{B_T}{\sqrt{N}}
    \end{split}
\end{equation}
where 
\be{CT2}
B_T =  \left[\l T + 1 + K_1 \sqrt{T(\l+r)}\right]e^{(\a+3\l + 2r)T},
\ee
and $K_1$ is the best constant in the Burkholder-Davis-Gundy inequality in $L^1$. Moreover, defining
\be{v2}
v_{N,2}(t) := \frac{1}{N} \sum_{i=1}^N x_i^2(t),
\ee
we have
\be{v2est}
\E\left[\sup_{t \in [0,T]} \left| v_{N,2}(t) - \E\left[\bar x^2(t)\right] \right| \right] \leq \frac{C_T}{\sqrt{N}},
\ee
where 
\be{CT3}
C_T = 2 A_T + (1 + \l T + K_1 (3 \sqrt{\l} + \sqrt{r})\sqrt{T})e^{(2\a+r)T}.
\ee

\end{corollary}

\section{Macroscopic dynamics of aggregated variables}\label{sec:macro_dyn}

\subsection{Evolution of the mean values}

By taking the expectation in Eq. \eqref{eq:prop:chaos}, it is easy to obtain the evolution equations for 
\begin{equation}
m(t) \coloneqq \P\left[\bar{x}(t) >0\right]
\end{equation}
and
\begin{equation}
v(t) \coloneqq \E\left[\bar{x}(t)\right], 
\end{equation}
which are given by
\begin{equation}\label{eq:m(t)v(t)}
\begin{split}
	& \dot{m}(t) = \l(1-m(t))v(t) - r m(t) \\
         & \dot{v}(t) = -\alpha v(t) + \l(1-m(t))v(t) - r v(t) .
\end{split}
\end{equation} 
It is easily seen that $(0,0)$ is always a fixed point for \eqref{eq:m(t)v(t)}. Under the condition $\frac{r+\alpha}{\l}<1$ a nonzero fixed point exists:
\begin{equation}\label{eq:mstar:vstar}
        \left(m^{\ast}, v^{\ast}\right) \coloneqq \left(1-\frac{r+\alpha}{\l}, \frac{r }{r+\alpha}\left(1-\frac{r+\alpha}{\l}\right) \right)    \end{equation} 
Linear stability for these fixed points is readily checked with the Jacobian matrix of the vector field in  \eqref{eq:m(t)v(t)}:
\begin{equation}
    J(m,v) = 
    \begin{bmatrix}
    -\l v-r & \l(1-m)\\
    -\l v & \l(1-m) - (r+\alpha)
    \end{bmatrix}
\end{equation}
As 
\begin{equation}
    J(0,0)=
    \begin{bmatrix}
    -r & \l\\
    0 & \l-(r+\alpha)
    \end{bmatrix}
\end{equation}
we see that the origin is linearly stable for $\frac{r+\alpha}{\l}>1$. Its asymptotic global stability for $\frac{r+\alpha}{\l} \geq 1$ follows from the absence of limit cycles, 
that can be proved by employing Dulac's criterion: defining $g(m,v)\coloneqq v^{-1}$, it holds
\begin{equation}
	div\left(g(m,v)\begin{bmatrix} \dot{m}\\ \dot{v}\end{bmatrix}\right) = -\l - 
	\frac{r}{v} < 0 \quad \forall\, 0\leq m\leq 1,\, 0 < v \leq 1.
\end{equation}
Under the condition $\frac{r+\alpha}{\l}<1$, the Jacobian matrix
\begin{equation}\label{eq:jacobian:mstar:vstar}
    J(m^{\ast},v^{\ast}) = 
    \begin{bmatrix}
    -\frac{r \l}{r+\alpha} & r+\alpha\\
    -\frac{r \l}{r+\alpha}\left(1-\frac{r+\alpha}{\l}\right) & 0
    \end{bmatrix}
\end{equation}
has eigenvalues
\begin{equation}
    \lambda_\pm = \frac{1}{2}\left(-\frac{r \l}{r+\alpha}\pm \frac{r \l}{r+\alpha} \sqrt{ 1-\frac{4 (r+\alpha)^2}{r \l}\left(1-\frac{r+\alpha}{\l}\right)}\right) 
\end{equation}
which have negative real part. Thus the fixed point $(m^{\ast},v^{\ast}) $ is linearly stable whenever it exists. Its global stability follows form the above Dulac's criterion. Note also that the eigenvalues $ \lambda_\pm$  are real and negative when
$$\frac{4 (r+\alpha)^2}{r \l }\left(1-\frac{r+\alpha}{\l} \right) \leq 1 $$
which happens when 
\begin{equation}\label{eq:cminus:cplus}
\l\leq \l_{-}\coloneqq \frac{2 (r+\alpha)^2}{r}\left(1-\sqrt{\frac{\alpha}{r+\alpha} }\right) \text{ or } \l\geq \l_+ \coloneqq \frac{2 (r+\alpha)^2}{r}\left(1+\sqrt{\frac{\alpha}{r+\alpha} }\right) 
\end{equation}
and are complex and conjugate with negative real part when
$$\l_- < \l < \l_+ $$
Notice that $r+\alpha < \l_{-}$.  Thus, $(m^{\ast},v^{\ast})$ is a stable node when $r+\alpha < \l\leq \l_-$ or $\l \geq \l_+$ and a stable spiral when $\l_- < \l < \l_+$.

\begin{remark} \label{rem:nonmark}
In the Contact process without dissipation there is correspondence between the Markovianity of the one dimensional process $m_N(t)$ and the fact that its limit is {\em causal}, i.e. solves a well posed one dimensional Cauchy problem. In the dissipative Contact process, the limit process $(m(t),v(t))$ is causal, but its pre-limit counterpart is not Markovian. 
Indeed, the presence of the term
\[
\int_{[0,t]\times [0,+\infty)} x_i(s^-) \pmb{1}_{[0, \ind_{\{x_i(s^-)>0\}} r]}(u) N^\downarrow_i(ds,du)
\]
in the evolution of $x_i$, produces  for $v_N$ a martingale with a quadratic variation containing a term proportional to 
\[
\frac{1}{N} \sum_{i=1}^N \int_0^t x_i^2(s) ds,
\]
which is not a function of $(m_N, v_N)$. The contribution of this term, however, vanishes in the limit as $N \ra +\infty$. It can be shown that the dissipative Contact process does not admit any finite dimensional reduction, in the sense that no finite dimensional statistics of the microscopic variables has the Markov property for all $N$.

\end{remark}

\subsection{The distribution of $\bar x$}

Our aim in this subsection is to determine the asymptotic distribution of the viral load. We denote by $(\rho^{\mu_0}_{\bar x(t)})_{t\geq 0}$ the flow of the distributions of $\bar x$, with the initial condition $\bar x(0) \sim \mu_0$. Moreover, let $(m(t),v(t))$ be the solution of \eqref{eq:m(t)v(t)} with initial conditions $m(0) = \P(\bar x(0) >0)$, $v(0) = \E[\bar x(0)]$.
\begin{theorem}\label{thm:limit:rho:lambda}
For any initial distribution $\mu_0$ on $[0,1]$, the distribution at time $t$ of the limit process $\bar x$  when the initial distribution is $\bar x(0)\sim \mu_0$ is given by 
\begin{equation}\label{eq:general:rho}
\begin{split}
\rho_{\bar x(t)}^{\mu_0}(d\bar x) & = \int_0^1 \rho^{\bar x_0}_{\bar x(t)}(d\bar x) \mu_0(d\bar x_0)
\end{split}
\end{equation}
with
\begin{equation}\label{eq:fundamental:rho:complete}
\begin{split}
& \rho^{\bar x_0}_{\bar x(t)}(d\bar x) \coloneqq  
e^{-r t}\pmb{1}_{\bar x_0>0} \delta_{\bar x_0 e^{-\alpha t}}(d\bar x)  + g_t(\bar x)\pmb{1}_{( e^{-\alpha t},1]}(\bar x) d\bar x + k(t) \delta_0(d\bar x)
\end{split}
\end{equation}
the marginal law at time $t$ of the limit $\bar x$ process started at the deterministic initial condition $\bar x(0)=\bar x_0\in [0,1]$, with 
\begin{equation}\label{eq:def:gtlam}
g_t(\bar x) \coloneqq \frac{\l \, k\left(t+\frac{1}{\alpha}\ln{\left(\bar x\right)}\right) v\left(t+\frac{1}{\alpha}\ln{\left(\bar x\right)}\right)}{\alpha }\bar x^{\frac{r-\alpha}{\alpha}}
\end{equation}
and $k(t)$ being the solution to 
\begin{equation}\label{eq:ode:k(t)}
\begin{cases}
\dot{k}(t) + \l v(t) k(t) = r\left(1-k(t)\right)\\
k(0) = \pmb{1}_{\{\bar x_0 = 0\}}
\end{cases}
\end{equation}
Besides $\delta_0$, the unique stationary distribution is given by
\begin{equation}\label{eq:lambda:stationary}
\begin{split}
\rho^{\ast}(d\bar x) & = \frac{ \l \, r}{\alpha } \frac{v^\ast}{\l \, v^\ast + r}\bar x^{\frac{r-\alpha}{\alpha}} d\bar x+ k^\ast \delta_0(d\bar x)\\
& = \frac{r}{\alpha} \left( 1- \frac{r+\alpha}{\l}\right) \bar x^{\frac{r-\alpha}{\alpha}} d\bar x + \frac{r+\alpha}{\l}\delta_0(d \bar x)
\end{split}
\end{equation}
where $k^{\ast} = \lim_{t\to +\infty}k (t) = \frac{r}{\l \, v^\ast + r}$ is the stationary solution to Eq. \eqref{eq:ode:k(t)}. Moreover, $\rho^{\ast}$ is the limit distribution (as $t \ra +\infty$) for all initial distributions $\mu_0$ different from $\delta_0$.
\end{theorem}

\begin{remark}\label{rmk:heuristic:rho}
We give a heuristic interpretation of \eqref{eq:fundamental:rho:complete}.  
The first summand arises as, starting from $ \bar x_0>0$, the limit process $ \bar x$ evolves deterministically and decreases exponentially fast to zero until a jump to $ \bar x=0$ occurs. This happens with rate $r$, so the probability that, at any time $t>0$ and starting from $ \bar x_0 >0$, the process $ \bar x$ has not jumped to zero yet is $e^{-r t}$. As soon as $ \bar x$ undergoes a jump to zero, the initial condition $ \bar x_0$ is forgotten.

This is when the second and the third summands come into play.  The second summand constitutes the absolutely continuous part of $\rho^{ \bar x_0}_{ \bar x(t)}$. It is turned on as $ \bar x$, after having reached zero for the first time, jumps for the first time to $1$. Indeed, the only way $\pmb{1}_{( e^{-\alpha t},1]}( \bar x(t))$ can be different from zero is when $ \bar x$ decays deterministically from $ \bar x_0$ for a time $s\leq t$, then jumps to zero and then to $1$ before time $t$, and from $1$ resumes a deterministic decay (possibly with other extra jumps to zero and $1$) so that, by time $t$, $ \bar x(t) \in ( e^{-\alpha t}, 1]$. In fact, if $ \bar x$ underwent only a deterministic decay from $ \bar x_0$ at time $0$ to time $t$, we could only find $ \bar x(t)\leq  e^{-\alpha t}$.    

The last summand accounts for the fact that the process can jump to zero and stay in zero for some time. 
\end{remark}

\section{The fluctuation process and its Gaussian limit}\label{sec:normal_fluct}

We now consider the fluctuation process
\begin{equation}\label{eq:def:fluctuations}
    \begin{split}
	X^N (t) \coloneqq 
	\begin{bmatrix}
         \xi^N(t)\\
	\eta^N(t) 
	\end{bmatrix} \coloneqq
	\begin{bmatrix}
	\sqrt{N}\left(m_N(t) - m(t)\right)\\
	\sqrt{N}\left(v_N(t) - v(t)\right)
	\end{bmatrix}
    \end{split}
\end{equation}
As observed in Remark \ref{rem:nonmark}, this process is not Markovian. However, we prove a Central Limit Theorem that implies its convergence to a Gauss-Markov process.

\begin{theorem}[\textbf{Diffusion approximation of the fluctuation process}]\label{thm:fluct:convergence}
As in Theorem \ref{thm:prop:chaos} we assume that $(x_i(0))_{i=1}^N$ are independent identically distributed random variables with values in $[0,1]$ and law $\mu_0$. In particular, $X^N(0)$ converges in distribution to a Gaussian random variable $X_0$.
Then $X^N$ converges in distribution, in any interval $[0,T]$, to the solution $X(t) = \left[\xi(t),\, \eta(t) \right]^T$ to the following linear stochastic differential equation
\begin{equation}\label{eq:fluct:limit}
    \begin{split}
        & d X(t) = 
        F (t)
        X(t) dt + 
        \Sigma(t) dW(t)\\
        & 
        X(0) = X_0
    \end{split}
\end{equation}
where  
\begin{equation}\label{eq:def:F}
F(t) \coloneqq 
\begin{bmatrix}
         -\left(\l v(t)+r\right) & \l(1-m(t))\\
        -\l v(t) & \l\left(1-m(t)\right) - r- \alpha
\end{bmatrix} 
\end{equation} 
is the drift matrix, $W(t)$ is a two-dimensional standard Brownian motion and the diffusion matrix $\Sigma(t)$ is the (unique, symmetric) square root of the matrix $A(t)$ given by
\begin{equation}\label{eq:def:A}
A(t) = \Sigma(t)\Sigma^T(t) = 
\begin{bmatrix}
\l \left(1-m(t)\right)v(t) +r m(t) & \l \left(1-m(t)\right)v(t) +r \l v(t) \\
\l\left(1-m(t)\right)v(t) +r v(t) & \l \left(1-m(t)\right) v(t) + r \E\left[\bar x^2(t)\right]
\end{bmatrix}
\end{equation}
where $\bar x(t)$ is the limit process defined in \eqref{eq:prop:chaos}.
\end{theorem}

\subsection{Power spectral density of the fluctuation process in the supercritical regime} \label{sec:specde}

The Gaussian process $X(t)$ defined in \eqref{eq:fluct:limit} describes, for $N$ large, the fluctuations of the process $(m_N, v_N)$ around its deterministic limit $(m,v)$. If we assume $\l > r+\a$ we have that, except for the trivial case $m(0) = v(0) = 0$,
\[
\lim_{t \ra +\infty} (m(t),v(t)) =    \left(m^{\ast}, v^{\ast}\right) \coloneqq \left(1-\frac{r+\alpha}{\l}, \frac{r}{r+\alpha}\left(1-\frac{r+\alpha}{\l}\right) \right)
\]
Moreover, the evolution of $X(t)$, as $t$ converges to infinity, converges to the stationary solution of the linear stochastic equation
\[
d X(t) = 
        F 
        X(t) dt + 
        \Sigma dW(t),
        \]
where 
\[
 F = 
 \lim_{t 
 \ra +\infty} F(t)  \coloneqq 
\begin{bmatrix}
         - \left(\l v^{\ast}+r\right) & \l\left(1-m^{\ast}\right) \\
        -\l v^{\ast} & 0
\end{bmatrix} 
\]
\[
\Sigma\Sigma^T = \lim_{t \ra +\infty} \Sigma(t)\Sigma^T(t) = 
\begin{bmatrix}
\l \left(1-m^\ast\right)v^\ast +r m^\ast & \l\left(1-m^\ast\right)v^\ast +r \l v^\ast \\
\l \left(1-m^\ast\right)v^\ast +r v^\ast & \l \left(1-m^\ast\right) v^\ast + r \E_{\rho^\ast}[ x^2]
\end{bmatrix}
\]
where the distribution $\rho^\ast$ is given in \eqref{eq:lambda:stationary} and $\E_{\rho^\ast}$ denotes the expectation taken with respect to it. 
By abuse of notation we still denote by $X(t)$ this stationary Gaussian process; it describes the long time fluctuations of the process $(m_N, v_N)$ around the equilibrium $ \left(m^{\ast}, v^{\ast}\right)$. Several features of this process can be obtained by its {\em power spectral density matrix} $S(\o)$ defined by
\begin{equation} \label{eq:specdendef}
S(\o) := 
\lim_{t \ra +\infty} \frac{1}{t} \E\left[ \hat{X}_t (\o)  \hat{X}^T_t (\o) \right],
\end{equation}
where, for $\o \in \R$, 
\[
 \hat{X}_t (\o) := \int_0^t X(s) e^{-i \o s} ds.
 \]
 In particular, the diagonal element $S_{11}(\o)$ (resp. $S_{22}(\o)$), provides the spectral profile of the frequencies of  the two components of $X(t)$.
 \begin{proposition} \label{prop:specden}
 Denoting by $(a_{ij})_{i,j =1,2}$ the entries of the matrix    $\Sigma\Sigma^T$ and by $(F_{ij})_{i,j =1,2}$ the entries of the matrix $F$, we have
\begin{equation}\label{eq:S_11(omega)}
    \begin{split}
        S_{11}(\omega)=\frac{a_{11} \omega^2 + a_{11} F^2_{22} + a_{22} F^{2}_{12} -2 F_{12}F_{22}a_{12}}{\left(\omega^2-det(F)\right)^2 + \omega^2 tr(F)^2}
    \end{split}
\end{equation}
\begin{equation}\label{eq:S_22(omega)}
    \begin{split}
        S_{22}(\omega)=\frac{a_{22} \omega^2 + a_{22} F^2_{11} + a_{11} F^{2}_{21} -2 F_{11}F_{21}a_{21}}{\left(\omega^2-det(F)\right)^2 + \omega^2 tr(F)^2}
    \end{split}
\end{equation}

\end{proposition}

The proof of this proposition can be found for instance in \cite{gardiner}, Section 4.4.6.  

\begin{figure} 
\centering
    \includegraphics[width=0.80\textwidth]{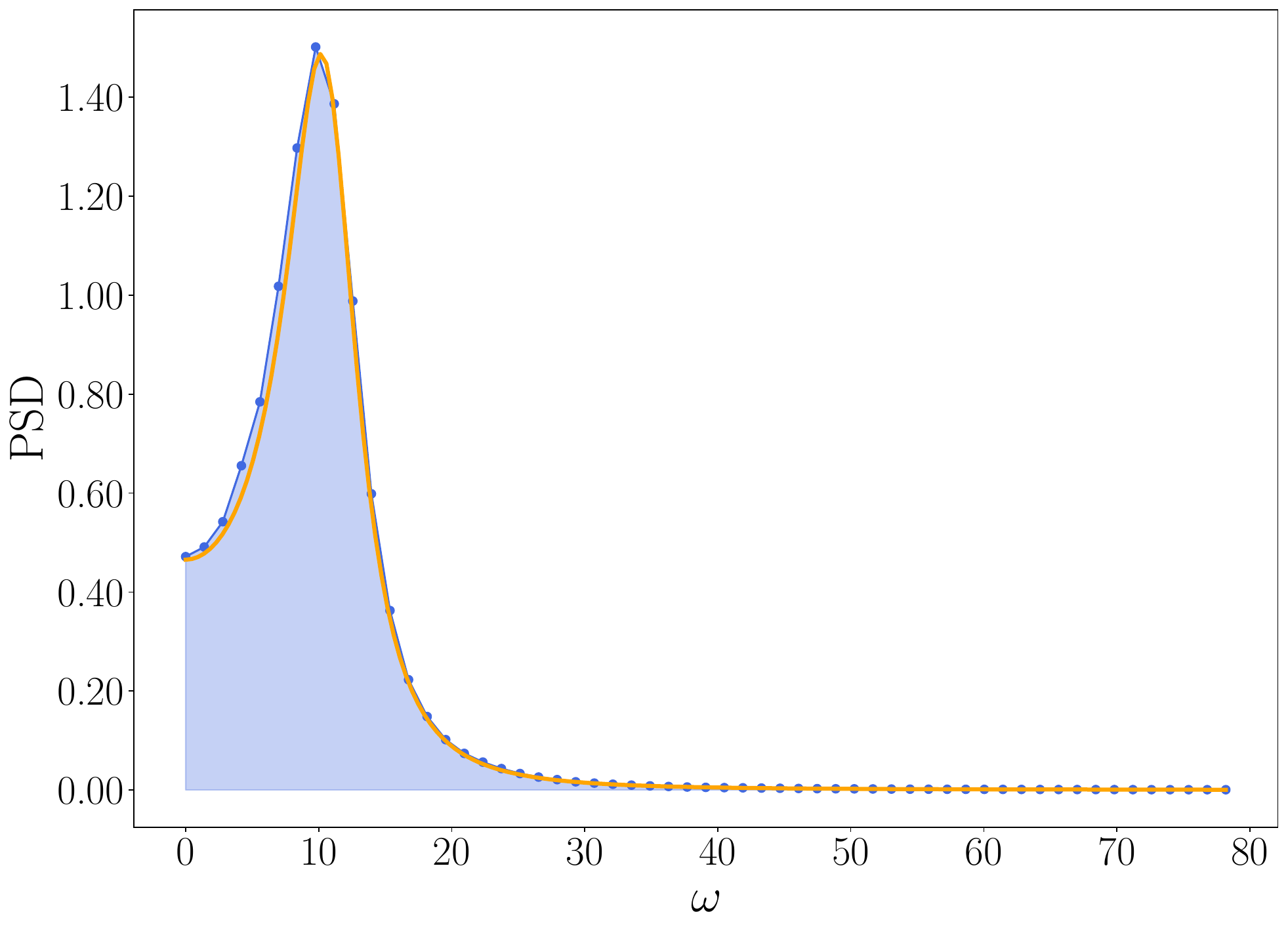}
     \caption{The spectral density of the process $\xi(t)$ in \eqref{eq:fluct:limit} (in orange) compared to the estimated spectrum of $\xi_N$ in \eqref{eq:def:fluctuations} (blue points) for the parameters $N=10000$, $r=5$, $\lambda=100$, $\rho= 0.7$. To obtain the blue curve, $100$ simulations of system \eqref{eq:micro:dyn} were performed employing an Euler scheme with time-step $0.001$.}
     \label{fig:fluct:spectrum}
\end{figure}

It is interesting to see the behavior of these spectral densities in the case the time scales $r^{-1}$ (the expected time needed for an infected individual to lose immunity) and $\a^{-1}$ (the time an infected individual remains contagious) are well separated, i.e. $\a \gg r$. To describe the behavior of the spectral density in this regime we keep $r$ and $\rho:= \frac{\a}{\lambda} \in (0,1)$ fixed and let $\lambda$, and therefore also $\a$, increase to infinity. Note that in this regime, using the notation in \eqref{eq:cminus:cplus}, $\lambda_- \sim \a = \rho \lambda < \lambda$ and $\lambda_+ \sim \frac{4 \rho^2 \lambda^2}{r} > \lambda$ for $\lambda$ large. So the condition $\lambda_- < \lambda < \lambda_+$ is satisfied, and the equilibrium $(m^{\ast}, v^{\ast})$ is a stable spiral. Moreover
\[
S_{11}(\omega) \sim \frac{2 r \left(1-\rho\right)\omega^2 + r\left(1-\rho\right)  \rho^2 \lambda^2}{\left(\omega^2 - r  \left(1-\rho\right)\lambda\right)^2 + \frac{r^2}{\rho^2}\omega^2 }
\]
and
\[
S_{22}(\omega) \sim \frac{r \left(1-\rho\right) \omega^2 + r^3\frac{1-\rho}{\rho^2}  + 2 r^3 \frac{\left(1-\rho\right)^3}{\rho^2}  -2 r^3 \frac{\left(1-\rho\right)^2}{\rho^2} }{\left(\omega^2 - r  \left(1-\rho\right)\lambda\right)^2 + \frac{r^2}{\rho^2}\omega^2 }.
\]
It is easily seen that, for $\lambda$ large, both spectral densities have a sharp peak around the frequency $\o^* := \sqrt{r(1-\rho)\lambda}$, as shown in Figure \ref{fig:fluct:spectrum}: random fluctuations around equilibrium exhibit oscillations with a nearly deterministic period. This analysis motivates the scaling limit performed in next section.

\section{Rescaling the parameters} \label{sec:rescale}

The aim of this section is to obtain the limit behavior of the microscopic dynamics  when we let the parameters $\a$ and $\l$ diverge with $N$. Consistently with what we have seen in Section \ref{sec:specde}, we fix $r>0$, $\rho\in (0,1)$, set $\a_N = \rho \l_N$ and let $\l_N \uparrow +
\infty$ as $N \uparrow +\infty$. We have observed that this asymptotics produces a sharp peak in the power spectrum. 
On the one hand, since the power spectral density reflects the long-time behavior of a process, we do not expect the peak of $S_{11}$ and $S_{22}$ to be due to the damped oscillations which the limit system \eqref{eq:m(t)v(t)} exhibits when far from equilibrium. On the other hand, for $\l \gg r$, by looking at the linearization of system \eqref{eq:m(t)v(t)} around its stable equilibrium we might expect that \eqref{eq:m(t)v(t)} behaves as a deterministic harmonic oscillator, the dampening factor acting on very long time scales. In turn, one might think that the sharp peak of the power spectrum essentially captures deterministic oscillations with a negligible dampening. However, numerical simulations of system \eqref{eq:def:fluctuations} show oscillations with a nearly deterministic period, but with a random amplitude, hence, different in nature from the ones of a deterministic oscillator. 
  To explain this phenomenon, taking into account that we are interested in the long time behavior, we begin by assuming the initial state $\left(x_i(0)\right)_{i = 1}^N$ to be a system of i.i.d. variables, as in Theorem \ref{thm:fluct:convergence}, but such that $\P(x_i(0) > 0 ) = m^{\ast}$, $\E[x_i(0)] = v^{\ast}$, i.e. the equilibrium values have been attained. For instance, we may assume $x_i(0)$ to be distributed according to the $\rho^{\ast}$ defined in \eqref{eq:lambda:stationary}. The key point now is to find a suitable rescaling of the fluctuations $\left(m_N(t) - m^{\ast}\right)$, $\left(v_N(t) - v^{\ast}\right)$ to have a nontrivial limit. Note that even propagation of chaos is not obvious anymore: indeed, since some parameters of the model diverge with $N$, the possibility that propagation of chaos deteriorates in short time cannot be ignored. Moreover, if oscillations are detected in the limit, a specific time rescaling is forced by the the previous observation that the asymptotic frequency is of order $\sqrt{\l}$.

\begin{theorem} \label{th:rescal}
Define the rescaled fluctuation processes as follows:
\begin{equation}\label{eq:def:fluctuations:rescaled}
    \begin{split}
        &\hat{\xi}^N (t) \coloneqq \sqrt{N}\frac{m_N\left(\frac{t}{\sqrt{\l_N}}\right) - m^\ast}{\l_N^{1/4} } \\
        & \hat{\eta}^N (t) \coloneqq \l_N^{1/4} \sqrt{N}\left(v_N\left(\frac{t}{\sqrt{\l_N}}\right) - v^{\ast} \right)    \end{split}
\end{equation}
and assume the initial state $\left(x_i(0)\right)_{i = 1}^N$ to be a system of i.i.d. variables such that $\P(x_i(0) > 0 ) = m^{\ast}$, $\E[x_i(0)] = v^{\ast}$. Moreover assume $\a_N = \rho \l_N$,  and
\[
\lim_{N \ra +\infty} \l_N = +\infty \hspace{2cm} \lim_{N \ra +\infty} \frac{\sqrt{\l_N}}{\log N} = 0.
\]
Then the process $(\hat{\xi}^N (t), \hat{\eta}^N (t))^T$ converges in distribution, in any interval $[0,T]$, to the random harmonic oscillator
\begin{equation} \label{eq:rescal}
\begin{split}
d\hat{\xi}(t) & = \rho \hat{\eta} (t) dt \\
d\hat{\eta} (t) & = -\frac{r}{\rho}(1-\rho) \hat{\xi}(t)dt + r(1-\rho)dW(t) \\
\hat{\xi}(0)  & = \hat{\eta}(0)  = 0
\end{split}
\end{equation}
where $W(t)$ is a standard Brownian motion.

\end{theorem}

\begin{figure} 
\centering
    \includegraphics[width=0.80\textwidth]{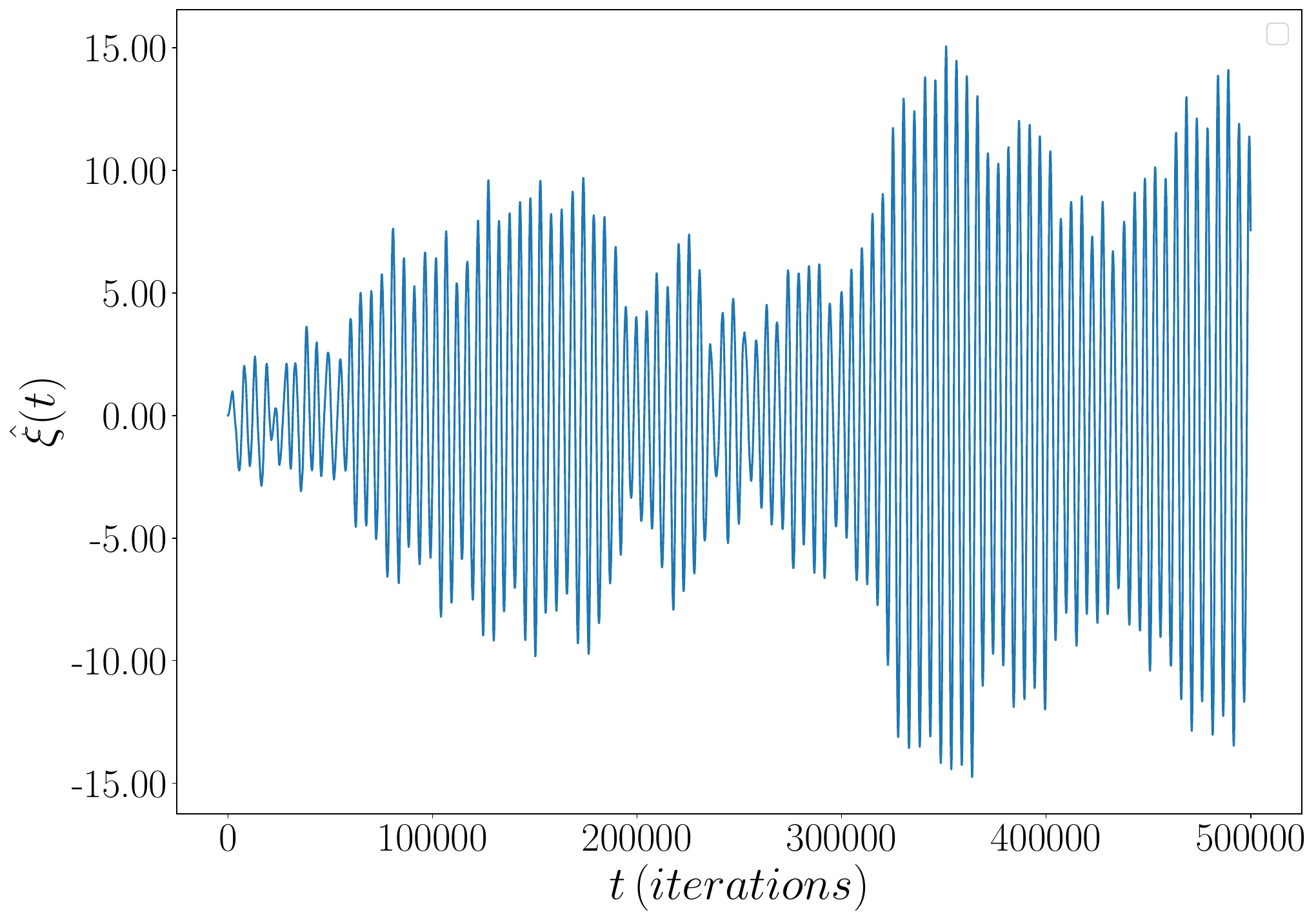}
     \caption{Sample trajectory of $\hat\xi(t)$ (system \eqref{eq:rescal}) at stationarity.}
     \label{fig:fluct:trajectories}
\end{figure}

As the simulation in Figure \ref{fig:fluct:trajectories} shows, the Brownian noise has little effects on the frequency of the oscillations, which is $\sqrt{r(1-\rho)}$ for the deterministic system, but substantially affects the amplitudes.

\section{Proofs}\label{sec:proofs}

\subsection{Proof of Theorem \ref{thm:prop:chaos}}

It is convenient to set $\s_i(t) := \ind_{\{x_i(t) > 0\}}$, so that \eqref{eq:micro:dyn} can  be split into the system
\begin{equation}\label{eq:micro:dyn2}
\begin{split}
     x_i(t) & = x_i(0) + \int_0^t -\alpha x_i(s) ds +  \int_{[0,t]\times [0,+\infty)} \pmb{1}_{[0, (1-\s_i(s^-))  \l v_N(s^-)]}(u) N^\uparrow_i(ds,du) \\
     &  - \int_{[0,t]\times [0,+\infty)} x_i(s^-) \pmb{1}_{[0,  \s_i(s^-) r ]}(u) N^\downarrow_i(ds,du) \\
     \s_i(t) & = \s_i(0) +  \int_{[0,t]\times [0,+\infty)} \pmb{1}_{[0, (1-\s_i(s^-))  \l v_N(s^-)]}(u) N^\uparrow_i(ds,du) \\
     &  - \int_{[0,t]\times [0,+\infty)} \pmb{1}_{[0, \s_i(s^-) r ]}(u) N^\downarrow_i(ds,du).
\end{split}
\end{equation}
A similar splitting applies to the limiting equation \eqref{eq:prop:chaos} if we define $\bar \s_i(t) := \ind_{\{\bar x_i(t) > 0\}} $:
\begin{equation} \label{reprlim}
  \begin{split}
         \bar x_i(t) & =  x_i(0) + \int_0^t -\alpha \bar x_i(s) ds +  \int_{[0,t] \times [0,+\infty)} \pmb{1}_{[0, (1-\bar\s_i(s^-))\l v(s^-)]}(u) N_i^\uparrow(ds,du) \\
         &  - \int_{[0,t] \times [0,+\infty)} \bar x_i(s^-) \pmb{1}_{[0, \bar\s_i(s^-) r ]}(u) N_i^\downarrow(ds,du) \\
         \bar\sigma_i(t) & = \sigma_i(0)  +  \int_{[0,t]\times [0,+\infty)} \pmb{1}_{[0,(1-\bar\sigma_i(s^-))\lambda v(s^-)]}(u)  N^\uparrow_i(ds,du)  \\ & - \int_{[0,t]\times [0,+\infty)} \pmb{1}_{[0,\bar\sigma_i(s^-) r ]}(u)  N^\downarrow_i(ds,du)
    \end{split}
 \end{equation}

It follows that
\be{PC1}
\begin{split}
| x_i(t) -  \bar x_i(t)| & \leq \alpha\int_0^t | x_i(s) -  \bar x_i(s)|ds + \int_{[0,t]\times [0,+\infty)} \left|  \pmb{1}_{[0, (1-\s_i(s^-))  \l v_N(s^-)]}(u) \right.\\ & \left.  - \pmb{1}_{[0,(1-\bar\s_i(s^-))\l v(s^-)]}(u)\right| N_i^\uparrow(ds,du) \\
& +  \int_{[0,t]\times [0,+\infty)} \left|x_i(s^-) \pmb{1}_{[0,  \s_i(s^-) r ]}(u)  - \bar x_i(s^-) \pmb{1}_{[0, \bar
\s_i(s^-)r ]}(u) \right| N_i^\downarrow(ds,du)
\end{split}
\ee
and similarly
\be{PC2}
\begin{split}
|\s_i(t) - \bar\sigma_i(t)| & \leq \int_{[0,t]\times [0,+\infty)} \left| \pmb{1}_{[0, (1-\s_i(s^-))  \l v_N(s^-)]}(u) - \pmb{1}_{[0,(1-\bar\sigma_i(s^-))\lambda v(s^-)]}(u) \right|  N^\uparrow_i(ds,du) \\
& + \int_{[0,t]\times [0,+\infty)}  \left|  \pmb{1}_{[0,  \s_i(s^-) r ]}(u) - \pmb{1}_{[0,\bar\sigma_i(s^-) r]}(u)   \right| N^\downarrow_i(ds,du).
\end{split}
\ee
Note that both in \eqref{PC1} and \eqref{PC2}, by increasingness of the integrals in the r.h.s,  the l.h.s. can be replaced by $\sup_{s \leq t} | x_i(s) -  \bar x_i(s)|$ and $\sup_{s \leq t} | \s_i(s) -  \bar \s_i(s)|$ respectively. Taking expectation and letting 
\[
\varphi(t) = \E\left[ \sup_{0 \leq s \leq t} \left( | x_i(s) -  \bar x_i(s)| + | \s_i(s) -  \bar \s_i(s)| \right) \right],\
\]
we obtain
\be{PC3}
\begin{split}
\varphi(t) & \leq \a \int_0^t \varphi(s) ds \\
& + \E \left[\int_{[0,t]\times [0,+\infty)} \left|  \pmb{1}_{[0, (1-\s_i(s^-))  \l v_N(s^-)]}(u)  - \pmb{1}_{[0,(1-\bar\s_i(s^-))\l v(s^-)]}(u)\right| N_i^\uparrow(ds,du) \right] \\
& + \E \left[ \int_{[0,t]\times [0,+\infty)} \left|x_i(s^-) \pmb{1}_{[0,  \s_i(s^-) r ]}(u)  - \bar x_i(s^-) \pmb{1}_{[0, \bar
\s_i(s^-) r ]}(u) \right| N_i^\downarrow(ds,du)\right] \\
& +  \E \left[\int_{[0,t]\times [0,+\infty)} \left| \pmb{1}_{[0, (1-\s_i(s^-))  \l v_N(s^-)]}(u) - \pmb{1}_{[0,(1-\bar\sigma_i(s^-))\lambda v(s^-)]}(u) \right|  N^\uparrow_i(ds,du) \right] \\
& + \E \left[\int_{[0,t]\times [0,+\infty)}  \left|  \pmb{1}_{[0,  \s_i(s^-) r ]}(u) - \pmb{1}_{[0,\bar\sigma_i(s^-) r]}(u)   \right| N^\downarrow_i(ds,du) \right]
\end{split}
\ee
We now use the following facts: if $f(s,u)$ is bounded, positive and predictable then 
\be{smoothing}
 \E \left[\int_{[0,t]\times [0,+\infty)} f(s,u) N_i^\uparrow(ds,du) \right] =  \E \left[ \int_{[0,t]\times [0,+\infty)}f(s,u) du ds \right]
 \ee
and, for $a,a',b,b' \geq 0$
\[
\int_0^{+\infty} \left| a \pmb{1}_{[0,b]}(u) - a' \pmb{1}_{[0,b']}(u) \right| du \leq a|b-b'| + b'|a-a'|.
\]
This yields
\be{PC4}
\begin{split}
\E &  \left[\int_{[0,t]\times [0,+\infty)} \left|  \pmb{1}_{[0, (1-\s_i(s^-))  \l v_N(s^-)]}(u)  - \pmb{1}_{[0,(1-\bar\s_i(s^-))\l v(s^-)]}(u)\right| N_i^\uparrow(ds,du) \right] \\
& = \l \E \left[ \int_0^t \left|(1-\s_i(s))v_N(s) - (1-\bar\s_i(s))v(s) \right|ds \right] \\  & \leq \l \E \left[ \int_0^t |\s_i(s) - \bar\s_i(s)|ds \right] + \l \E \left[ \int_0^t |v_N(s) - v(s)| ds \right] \\
& \leq \l \int_0^t \varphi(s) ds + \l \E \left[ \int_0^t |v_N(s) - v(s)| ds \right]
\end{split} 
\ee
Similarly 
\be{PC5}
\begin{split}
\E & \left[ \int_{[0,t]\times [0,+\infty)} \left|x_i(s^-) \pmb{1}_{[0,  \s_i(s^-) r ]}(u)  - \bar x_i(s^-) \pmb{1}_{[0, \bar
\s_i(s^-) r ]}(u) \right| N_i^\downarrow(ds,du)\right] \\ & \leq r \E\left[ \int_0^t |x_i(s) - \bar x_i(s)| ds \right] + r \E\left[ \int_0^t |\s_i(s) - \bar \s_i(s)| ds \right] \\ & \leq r \int_0^t \varphi(s) ds,
\end{split}
\ee
Estimating in the same way the other two terms in the r.h.s. of \eqref{PC3}, from \eqref{PC3}, \eqref{PC4} and \eqref{PC5} we obtain
\be{PC6}
\varphi(t) \leq (\a+2\l + 2r) \int_0^t \varphi(s) ds + \l \E \left[ \int_0^t |v_N(s) - v(s)| ds \right].
\ee
Noticing that, letting $\bar v_N(s) \coloneqq \frac{1}{N}\sum_{i=1}^N \bar x_i(s)$, 
\begin{equation*}
\begin{split}
\E \left[|v_N(s) - v(s)|\right] & \leq \E\left[| v_N(s) - \bar v_N(s)|\right] + \E\left[|\bar v_N(s) - v(s)| \right] \leq \phi(s) + \left( \mbox{Var}(\bar v_N(s)) \right)^{1/2} \\
				       & \leq   \phi(s) + \frac{1}{\sqrt{N}}, 
\end{split}
\end{equation*}
we finally obtain from \eqref{PC6}  
\[
\varphi(t) \leq (\a+3\l + 2r) \int_0^t \varphi(s) ds + \frac{\l T }{\sqrt{N}}.
\]
A direct application of Gromwall Lemma provides
\[
\varphi(T) \leq \frac{\l T }{\sqrt{N}} e^{(\a+3\l + 2r)T},
\]
and the proof is complete.

\subsection{Proof of Corollary \ref{cor:LLN}}

By adding and subtracting $\bar m_N(t) \coloneqq \frac{1}{N}\sum_{i=1}^N \bar\sigma_i(t)$, 
\begin{equation}\label{eq:cor:LLN:step1}
    \begin{split}
        \E\left[\sup_{t\in[0,T]} \left| m_N(t) - \P[\bar x(t) >0] \right|\right] & \leq  \E\left[\sup_{t\in[0,T]} \left| m_N(t) - \bar m_N(t)\right|\right]+\E\left[\sup_{t\in[0,T]} \left|\bar m_N(t) - m(t)\right|\right] \\
        & \leq \frac{1}{N}\sum_{i=1}^N\E\left[\sup_{t\in [0,T]}  \left|\sigma_i(t) -\bar\sigma_i(t)\right|\right]  +\E\left[\sup_{t\in[0,T]} \left|\bar m_N(t) - m(t)\right|\right]
    \end{split}
\end{equation}
The first term in the r.h.s. of \eqref{eq:cor:LLN:step1}, by Theorem \ref{thm:prop:chaos}, is bounded by $\frac{A_T}{\sqrt{N}}$. So we show that a similar bound holds for the second term in the r.h.s. of \eqref{eq:cor:LLN:step1}. Averaging over $i$ in \eqref{reprlim} we get
\begin{equation}\label{221}
	\bar m_N(t) = \bar m_N(0) + \int_0^t \left( \lambda v(s) (1-\bar m_N(s)) - r \bar m_N(s) \right) ds + \frac{1}{N} M_N(t)
\end{equation}
with
\be{martMN}
M_N(t) = \sum_{i=1}^N \left[\int_{[0,t]\times [0,+\infty)} \pmb{1}_{[0,(1-\bar\sigma_i(s^-))\lambda v(s^-)]}(u)  \tilde N^\uparrow_i(ds,du) - \int_{[0,t]\times [0,+\infty)} \pmb{1}_{[0,\bar\sigma_i(s^-) r ]}(u) \tilde N^\downarrow_i(ds,du)\right],
\ee
where we have introduced the compensated Poisson random measures defined by 
\begin{equation} \label{compPRM}
\begin{split}
\tilde{N}^\uparrow_i(ds,du) & := N^\uparrow_i(ds,du) - dsdu \\
\tilde{N}^\downarrow_i(ds,du) & := N^\downarrow_i(ds,du) - dsdu. 
\end{split}
\end{equation}
We recall that for a predictable, positive and bounded $f(s,u)$ such that  
\[
\E\left[\int_{[0,T] \times [0,+\infty)} f^2(s,u)dsdu \right]<+\infty,
\]
the integrals 
\[
I_i(t) := \int_{[0,t] \times [0,+\infty)} f(s,u) \tilde N^\uparrow_i(ds,du)
\]
 define orthogonal martingales with quadratic variation $ \int_{[0,t] \times [0,+\infty)} f^2(s,u)dsdu$.
It follows that $M_N(t)$ is a martingale having quadratic variation equal to $\sum_{i=1}^N \int_0^t [\lambda v(s) (1-\bar\sigma_i(s)) + r \bar\sigma_i(s)] ds$. From \eqref{eq:m(t)v(t)} we also obtain
\begin{equation}\label{222}
	m(t) = m(0) + \int_0^t \left( \lambda v(s) (1-m(s)) - r m(s) \right) ds.
\end{equation}
From \eqref{221} and \eqref{222} we get
\begin{equation} \label{223}
	\begin{split}
		| \bar m_N(t) - m(t) | & \leq  |\bar m_N(0) - m(0)| + \int_0^t ( \lambda v(s) + r) |\bar m_N(s) - m(s)|  ds + \frac{1}{N} |M_N(t)| \\ & \leq | \bar m_N(0) - m(0)| + \int_0^t ( \lambda + r) |\bar m_N(s) - m(s)|  ds + \frac{1}{N} |M_N(t)|.
	\end{split}
\end{equation}
Observe that, by independence of the components,
\[
\E\left[ |\bar m_N(0) - m(0)| \right] \leq \left(\mbox{Var}(\bar m_N(0)) \right)^{1/2} \leq \frac{1}{\sqrt{N}}.
\]
Letting 
\[
\varphi(t) := \E\left[ \sup_{s \in [0,t]} | \bar m_N(s) - m(s) | \right],
\]
using the Burkholder-Davis-Gundy inequality in $L^1$ and Jensen's inequality we obtain from \eqref{223}
\[
\begin{split}
\varphi(t) & \leq  \frac{1}{\sqrt{N}} + ( \lambda + r) \int_0^t \varphi(s)ds +  \frac{1}{N}\E\left[ \sup_{t\in [0,T]} |M_N (t)|\right] \\ 
&  \leq \frac{1}{\sqrt{N}} +( \lambda + r) \int_0^t \varphi(s)ds +  \frac{K_1}{N}\E\left[\left(\sum_{i=1}^N \int_0^t [\lambda v(s) (1-\bar\sigma_i(s)) + r \bar\sigma_i(s)] ds \right)^{1/2}\right]  \\
& \leq \frac{1}{\sqrt{N}} +( \lambda + r) \int_0^t \varphi(s)ds + \frac{K_1}{\sqrt{N}} \left( \int_0^T [\lambda v(s) (1- \bar m_N(s)) + r \bar m_N(s)] ds\right)^{1/2} \\
& \leq \frac{1}{\sqrt{N}} +( \lambda + r) \int_0^t \varphi(s)ds + \frac{K_1 \sqrt{T(\l+r)}}{\sqrt{N}}.
\end{split}
\]
giving
\be{estmbar}
\varphi(T)\leq \frac{1+K_1 \sqrt{T(\l+r)}}{\sqrt{N}} e^{( \lambda + r) T}.
\ee
By this estimate, \eqref{eq:cor:LLN:step1} and Theorem \ref{thm:prop:chaos}, we obtain the desired estimate for $ \E\left[\sup_{t\in[0,T]} \left| m_N(t) - \P\left[\bar{x}(t) >0\right]\right|\right]$. The estimate for $\E\left[\sup_{t \in [0,T]} \left| v_N(t) - \E\left[\bar{x}(t)\right] \right| \right] $ follows exactly the same steps, and it is omitted. So we are left to show \eqref{v2est}. We begin by observing that, letting $\bar v_{N,2}(t) \coloneqq \frac{1}{N}\sum_{i=1}^N \bar x^2_i(t)$, 
\be{v2split}
\begin{split}
\E\left[\sup_{t \in [0,T]} \left| v_{N,2}(t) - \E\left[\bar x^2(t)\right] \right| \right]  & \leq \frac{1}{N} \sum_{i=1}^N \E\left[\sup_{t \in [0,T]} |x_i^2(t) - \bar x_i^2(t)| \right] + \E\left[\sup_{t \in [0,T]} \left| \bar v_{N,2}(t) - \E\left[\bar x^2(t)\right] \right| \right] \\ & \leq  \frac{2}{N} \sum_{i=1}^N \E\left[\sup_{t \in [0,T]} |x_i(t) - \bar x_i(t)| \right]  + \E\left[\sup_{t \in [0,T]} \left| \bar v_{N,2}(t) - \E\left[\bar x^2(t)\right] \right| \right] \\
& \leq 2\frac{A_T}{\sqrt{N}} + \E\left[\sup_{t \in [0,T]} \left| \bar v_{N,2}(t) - \E\left[\bar x^2(t)\right] \right| \right] ,
\end{split}
\ee
where we have used estimate \eqref{eq:bound:prop:chaos}. To estimate $\E\left[\sup_{t \in [0,T]} \left| \bar v_{N,2}(t) - \E\left[\bar x^2(t)\right] \right| \right] $
we first apply It{\^o}'s formula to \eqref{eq:prop:chaos} and obtain
\begin{equation}\label{eq:ito:bar_x^2}
	\begin{split}
		\bar x^2_i(t) & = \bar x^2_i(0) -\alpha \int_0^t 2 \bar x^2_i(s) ds + \int_{[0,t]\times [0,+\infty)} \left[(\bar x_i(s^-) + 1)^2 - \bar x^2_i(s^-)\right] \pmb{1}_{[0,(1-\bar\sigma_i(s^-))\lambda v(s^-)]}(u) N^\uparrow_i(ds,du)  \\
		& + \int_{[0,t]\times [0,+\infty)}  - \bar x^2_i(s^-) \pmb{1}_{[0,\bar\sigma_i(s^-) r]}(u) N^\downarrow_i(ds,du) \\
		& = \bar x^2_i(0) + \int_0^t\left[ - 2 \alpha \bar x^2_i(s) + \lambda (1-\bar \sigma_i(s)) v(s) - r \bar x^2_i(s) \right] ds \\
		& + \int_{[0,t]\times [0,+\infty)} \left(2 \bar x_i(s^-) + 1\right) \pmb{1}_{[0,(1-\bar\sigma_i(s^-))\lambda v(s^-)]}(u) \tilde N^\uparrow_i(ds,du)  \\
		& - \int_{[0,t]\times [0,+\infty)}   \bar x^2_i(s^-) \pmb{1}_{[0,\bar\sigma_i(s^-) r]}(u) \tilde N^\downarrow_i(ds,du) .
	\end{split}
\end{equation}
On the one hand, averaging Eq. \eqref{eq:ito:bar_x^2} over $i$, we obtain $\bar v_{N,2}(t)$. 
On the other hand, taking the expectation in Eq. \eqref{eq:ito:bar_x^2}, we obtain the equation for $ \E[\bar x^2(t)]$: 
\begin{equation} \label{x2bar}
	\begin{split}
		\E[\bar x^2(t)] & = \E[\bar x^2(0)] -(2\alpha +r ) \int_0^t\E[\bar x^2(s)] ds + \int_0^t \lambda (1- m(s)) v(s) ds
	\end{split}
\end{equation}
Hence, defining
\[
\varphi(t) := \E\left[\sup_{s \in [0,t]} | \bar v_{N,2}(s) - \E[\bar x^2(s)] | \right],
\]
we obtain
\begin{equation} \label{v2semimart}
	\begin{split}
		\varphi(t) & \leq \E[|\bar v_{N,2}(0) - \E[\bar x^2(0)] | ] + (2\alpha + r) \int_0^t \varphi(s) ds + \lambda \int_0^t \E\left[ | m(s) - \bar m_N(s)|\right] ds\\
		& + \E\left[\sup_{s' \in [0,t]}\left|\frac{1}{N}\sum_{i=1}^N \int_{[0,s']\times [0,+\infty)} \left(2 \bar x_i(s^-) + 1\right) \pmb{1}_{[0,(1-\bar\sigma_i(s^-))\lambda v(s^-)]}(u) \tilde N^\uparrow_i(ds,du) \right|\right] \\
		& + \E\left[\sup_{s' \in [0,t]}\left|\frac{1}{N}\sum_{i=1}^N \int_{[0,s']\times [0,+\infty)}  \bar x^2_i(s^-) \pmb{1}_{[0,\bar\sigma_i(s^-) r]}(u) \tilde N^\downarrow_i(ds,du) \right|\right] 
\end{split}
\end{equation}
First we note that
\[
\int_0^t \E\left[ | m(s) - \bar m_N(s)|\right] ds \leq \int_0^t \left(\mbox{Var}(\bar m_N(s))\right)^{1/2}ds \leq \frac{T}{\sqrt{N}}.
\]
The two martingale terms in \eqref{v2semimart} can be bounded by the Burkholder-Davis-Gundy inequality:
\begin{multline*}
E\left[\sup_{s' \in [0,t]}\left|\frac{1}{N}\sum_{i=1}^N \int_{[0,s']\times [0,+\infty)} \left(2 \bar x_i(s^-) + 1\right) \pmb{1}_{[0,(1-\bar\sigma_i(s^-))\lambda v(s^-)]}(u) \tilde N^\uparrow_i(ds,du) \right|\right] \\
\leq K_1 \E\left[\left(\int_0^t \frac{1}{N^2}\sum_{i=1}^N \left(2 \bar x_i(s) + 1\right)^2 (1-\bar\sigma_i(s))\lambda v(s) ds \right)^{1/2}\right]  \leq \frac{3  K_1 \sqrt{\l T}}{\sqrt{N}}.
\end{multline*}
\begin{multline*}
 \E\left[\sup_{s' \in [0,t]}\left|\frac{1}{N}\sum_{i=1}^N \int_{[0,s']\times [0,+\infty)}  \bar x^2_i(s^-) \pmb{1}_{[0,\bar\sigma_i(s^-) r]}(u) \tilde N^\downarrow_i(ds,du) \right|\right]  \leq K_1 \E\left[\left(\int_0^T \frac{1}{N^2}\sum_{i=1}^N \bar x^4_i(s) r ds \right)^{1/2}\right]  \\ \leq \frac{K_1 \sqrt{r T}}{\sqrt{N}}.
 \end{multline*}
Thus we obtain
 \[
\varphi(t) \leq (2\alpha + r) \int_0^t \varphi(s)ds + \frac{D_T}{\sqrt{N}},
 \]
 where
 \[
 D_T = 1+ \l T + K_1 (3 \sqrt{\l} + \sqrt{r})\sqrt{T}
 \]
 
 from which
 \begin{equation} \label{C_Test2}
	\E\left[\sup_{t \in [0,T]} | \bar v_{N,2}(t) - \E[\bar x^2(t)] |\right] \leq \frac{D_T}{\sqrt{N}} e^{(2\a+r)T}
\end{equation}
follows by Gromwall Lemma. Going back to \eqref{v2split} we have
\[
\E\left[\sup_{t \in [0,T]} \left| v_{N,2}(t) - \E\left[\bar x^2(t)\right] \right| \right]  \leq \frac{2 A_T + D_T  e^{(2\a+r)T}}{\sqrt{N}},
\]
which completes the proof.

\subsection{Proof of Theorem \ref{thm:limit:rho:lambda}}
Recall that the limit process $\bar x = (\bar x(t))_{t\geq 0}$ is a piecewise deterministic Markov process.  
We first show that, for every $\bar x_0 \in [0,1]$ and $t\geq 0$, $\rho^{\bar x_0}_{\bar x(t)}$ is indeed a probability. Non negativity of $k(t)$ follows from \eqref{eq:ode:k(t)} observing that $\dot k(t)$ would be strictly positive if $k(t) < 0$. Thus we only have to check that 
\begin{equation} \label{eq:intmu}
\int_{e^{-\a t}}^1 g_t(x)dx = 1 - k(t) - e^{-r t}\pmb{1}_{\bar x_0>0}.
\end{equation}
Using the fact that, from \eqref{eq:ode:k(t)}, 
\[
\l e^{rt} v(t) k(t) = \frac{d}{dt}(1-k(t))e^{rt}
\]
we have, by the change of variable $y = t + \frac{1}{\a} \ln x$:
\[
\int_{e^{-\a t}}^1 g_t(x)dx = e^{-rt} \int_0^t \l e^{ry} k(y) v(y) dy = 1 - k(t) - e^{-rt}(1-k(0)),
\]
from which \eqref{eq:intmu} follows.

Now, it suffices to show that, for every $\bar x_0\in [0,1]$, the distribution $\rho^{\bar x_0}_{\bar x(t)}$ given in \eqref{eq:fundamental:rho:complete} is such that for every differentiable $f:[0,1] \ra \R$
\begin{equation}\label{eq:fundamental:eq}
\int_0^1 \mathcal{L}_t f(\bar x) \rho^{\bar x_0}_{\bar x(t)}(d\bar x) = \frac{d}{dt} \int_0^1 f(\bar x)\rho^{\bar x_0}_{\bar x(t)}(d\bar x) 
\end{equation}
where
\begin{equation}\label{eq:generator:lambdabar}
    \begin{split}
        \mathcal{L}_t f(\bar x) = -\alpha \bar x f'(\bar x) + \lambda v(t)\left[ f(1) - f(0)\right]\pmb{1}_{\{\bar x=0\}} + r \left[ f(0) - f(\bar x)\right]\pmb{1}_{\{\bar x >0\}}
    \end{split}
\end{equation}
 is the (time-dependent) infinitesimal generator of the process $\bar x$. We show the details for $\bar x_0 > 0$; the case $\bar x_0 = 0$ is similar and it is omitted.
 
 Employing Eq.s \eqref{eq:generator:lambdabar} and \eqref{eq:fundamental:rho:complete} we have that 
\begin{equation}\label{eq:lhsKolmog}
\begin{split}
\int_0^1 \mathcal{L}_tf(\bar x) \rho^{\bar x_{0}}_{\bar x(t)}(d\bar x) & = \int_0^1 \left(-\alpha\bar x f'(\bar x)  + \lambda v(t)\left[f(1)-f(0)\right] \pmb{1}_{\{\bar x=0\}} + r \left[f(0)-f(\bar x)\right]\pmb{1}_{\{\bar x >0\}} \right) \rho^{\bar x_{0}}_{\bar x(t)}(d\bar x) \\
& = -\alpha e^{-r t} \bar x_{0} e^{-\alpha t}f'(\bar x_{0} e^{-\alpha t}) -\alpha \int_{ e^{-\alpha t}}^1 \bar x f'(\bar x) g_t(\bar x) d\bar x \\
& + \lambda v(t) \left[f(1)-f(0)\right] e^{-r t}\int_0^1 \pmb{1}_{\{\bar x=0\}} \delta_{\bar x_{0} e^{-\alpha t}}(d\bar x) + \lambda v(t) \left[f(1)-f(0)\right] k(t) \\
& + r e^{-r t} \int_0^1 \left[f(0)-f(\bar x)\right]\pmb{1}_{\{\bar x >0\}}\delta_{\bar x_{0} e^{-\alpha t}}(d\bar x) + r \int_{ e^{-\alpha t}}^1 \left[f(0)-f(\bar x)\right] \pmb{1}_{\{\bar x >0\}} g_t(\bar x)d\bar x \\
& = -\alpha \bar x_{0} e^{-(\alpha+r) t}f'(\bar x_{0} e^{-\alpha t}) -\alpha \left(\left[\bar x f(\bar x) g_t(\bar x)\right]_{ e^{-\alpha t}}^1 - \int_{e^{-\alpha t}}^1 f(\bar x) \left(g_t(\bar x) + \bar x \partial_{\bar x} g_t(\bar x) \right)d\bar x \right)\\
& +  \lambda v(t) \left[f(1)-f(0)\right] k(t)  \\
& + r e^{-r t} \left[f(0)-f(\bar x_0 e^{-\alpha t})\right] + r \int_{e^{-\alpha t}}^1 \left[f(0)-f(\bar x)\right]  g_t(\bar x)d\bar x 
\end{split}
\end{equation}
On the other hand, we have that
\begin{equation}\label{eq:rhsKolmog}
\begin{split}
\frac{d}{dt}\int_0^1 f(\bar x) \rho^{\bar x_0}_{\bar x(t)}(d\bar x) & = \frac{d}{dt}\int_0^1 f(\bar x) e^{- r t}\delta_{\bar x_0 e^{-\alpha t}}(d\bar x) +\frac{d}{dt} \int_{e^{-\alpha t}}^1 f(\bar x)  g_t(\bar x)d\bar x + \dot k(t) f(0)\\
& = \frac{d}{dt}\left[e^{-r t} f(\bar x_0 e^{-\alpha t})\right] + \int_{e^{-\alpha t}}^1 f(\bar x) \partial_t g_t(\bar x)d\bar x +\alpha  e^{-\alpha t} f( e^{-\alpha t}) g_t(e^{-\alpha t})\\
&  + \dot k(t) f(0)\\
& = -r e^{-r t }f(\bar x_0 e^{-\alpha t})  -\alpha \bar x_0 e^{-(\alpha+r) t} f'(\bar x_0 e^{-\alpha t}) + \int_{e^{-\alpha t}}^1 f(\bar x) \partial_t g_t(\bar x)d\bar x \\
&+\alpha  e^{-\alpha t} f( e^{-\alpha t}) g_t(e^{-\alpha t}) + \dot k(t) f(0)
\end{split}
\end{equation}
Thus
\begin{equation}\label{eq:compare}
\begin{split}
\frac{d}{dt}\int_0^1 f(\bar x) \rho^{\bar x_0}_{\bar x(t)}(d\bar x) - \int_0^1 \mathcal{L}_tf(\bar x) \rho^{\bar x_{0}}_{\bar x(t)}(d\bar x) & = \int_{e^{-\alpha t}}^1 f(\bar x) \left[ \partial_t g_t(\bar x)  - (\a-r) g_t(\bar x) - \a \bar x \partial_{\bar x} g_t(\bar x) \right]d\bar x \\
& + f(0) \left[ \dot k(t) - re^{-rt} - r \int_{e^{-\a t}}^1 g_t(\bar x)d\bar x + \l v(t) k(t)\right] \\ & + f(1) \left[ \a g_t(1) - \l v(t) k(t) \right]  \\
& = \int_{e^{-\alpha t}}^1 f(\bar x) \left[ \partial_t g_t(\bar x) - (\a-r) g_t(\bar x) - \a \bar x \partial_{\bar x} g_t(\bar x) \right]d\bar x \\
& + f(0) \left[ \dot k(t) - re^{-rt} - r (1-k(t) - e^{-rt}) + \l v(t) k(t)\right] \\ & + f(1) \left[ \a g_t(1) - \l v(t) k(t) \right] \\
& = \int_{e^{-\alpha t}}^1 f(\bar x) \left[ \partial_t g_t(\bar x) - (\a-r) g_t(\bar x) - \a \bar x \partial_{\bar x} g_t(\bar x) \right]d\bar x \\
& + f(0) \left[ \dot k(t) - r (1-k(t)) + \l v(t) k(t)\right] \\ & + f(1) \left[ \a g_t(1) - \l v(t) k(t) \right] 
\end{split}
\end{equation}
By assumption $ \dot k(t) - r (1-k(t)) + \l v(t) k(t) = 0$. Moreover the equation
\[
\begin{split}
\partial_t g_t(\bar x) & =  (\a-r) g_t(\bar x) + \a \bar x \partial_{\bar x} g_t(\bar x) \\
g_t(1) & = \frac{\l}{\a} v(t) k(t)
\end{split}
\]
can be solved by the method of characteristics, and it is easy to checked that it is solved by \eqref{eq:def:gtlam}. This completes the proof of \eqref{eq:fundamental:eq}. The remaining statements concerning the limiting distribution are simple consequences of the stability of the fixed point $k^*$ of \eqref{eq:ode:k(t)}.

\subsection{Proof of Theorem \ref{thm:fluct:convergence}}

We will employ the following theorem: 
\begin{theorem}[\textbf{Diffusion approximation (Theorem VII, 4.1 \cite{EtKu2005})}]\label{thm:diffusion:approx}
Let $\left(X^N\right)_{N\in \mathbb{N}}$ and  $\left(B^N\right)_{N\in \mathbb{N}}$ $\R^d$-valued processes with càdlàg sample paths and let $A^N = \left(A^N_{i j}\right)$ be a symmetric $d\times d$ matrix-valued process such that $A^N_{i j}$ has càdlàg sample paths in $\R$ and $A^N(t)-A^N(s)$ is non-negative definite for all $t> s \geq 0$. 
Let $\mathcal{F}^N_t \coloneqq \sigma\left(X^N(s), B^N(s), A^N(s)\, : \, s\leq t \right)$. 
Let $\tau^N_h \coloneqq \inf\left\{t > 0\, : \, \vert X^N(t)\vert\geq h \, \text{or}\, \vert X^N(t^-)\vert \geq h \right\}$. 
Assume that
\begin{itemize}
    \item $M^N \coloneqq X^N - B^N$ and $M^N_i M^N_j - A^N_{i j }$ $i,j=1,\dots,d$ are $\left(\mathcal{F}^N_t\right)$-local martingales
    \item for each $T>0$, $h>0$
    \begin{equation}\label{eq:hyp4.4}
        \lim_{N \to +\infty}\E\left[ \sup_{t \leq T \wedge \tau^N_h} \vert B^N(t) - B^N(t^-)\vert^2\right] = 0
    \end{equation}
    \item for each $T>0$, $h>0$, $i,j=1,\dots,d$
    \begin{equation}\label{eq:hyp4.5}
        \lim_{N \to +\infty}\E\left[ \sup_{t \leq T \wedge \tau^N_h} \vert A^N_{i j}(t) - A^N_{i j}(t^-)\vert\right] = 0 
    \end{equation}
    \item for each $T>0$, $h>0$ 
    \begin{equation}\label{eq:hyp4.3}
        \lim_{N \to +\infty}\E\left[ \sup_{t \leq T \wedge \tau^N_h} \vert X^N(t) - X^N(t^-)\vert^2\right] = 0 
    \end{equation}
    \item there exist a continuous, symmetric, non-negative definite $d\times d$ matrix-valued function on $\R^d$, $a=\left(a_{i j}\right)$, and a continuous function $b: \R^d \to \R^d$ such that, for each $h>0$, $T>0$ and $i,j=1,\dots,d$, and for all $\epsilon>0$, 
    \begin{equation}\label{eq:hyp4.7}
        \lim_{N\to +\infty}\mathbb{P}\left(\sup_{t\leq T\wedge \tau^N_h} \Bigg\vert A^N_{ i j}(t) - \int_0^t a_{i j }(X^N(s))ds \Bigg\vert > \epsilon \right) = 0 
    \end{equation}
    and 
    \begin{equation}\label{eq:hyp4.6}
        \lim_{N\to +\infty}\mathbb{P}\left(\sup_{t\leq T \wedge \tau^N_h} \Bigg\vert B^N_i(t) - \int_0^t b_i(X^N(s))ds \Bigg\vert > \epsilon \right) = 0 
    \end{equation}
    \item the $\mathcal{C}_{\R^d}\left([0,+\infty)\right)$ martingale problem for
	\begin{equation}\label{A:EtKu}
	\tilde{A} \coloneqq \left\{\left(f,Gf\coloneqq \frac{1}{2}\sum_{i,j} a_{i j}\partial_i \partial_j f + \sum_i b_i \partial_i f\right)\, : \, f \in \mathcal{C}^{\infty}_c\left(\R^d\right) \right\} 
	\end{equation}
is well-posed. 

    \item the sequence of the initial laws of the $X^N$s converges in distribution to some probability distribution on $\R^d$, $\nu$. 
\end{itemize}
Then $\left(X^N\right)_N$ converges in distribution to the solution of the martingale problem for $(\tilde{A},\nu)$. 
That is, the laws of the processes $\left(X^N\right)_N$ converge weakly to the law of a process $X$ which is a weak solution of the SDE 
\begin{equation}
    \begin{split}
        dX(t) = b(X(t))dt +\Sigma(X(t)) dW(t)
    \end{split}
\end{equation}
where $b = (b_i)_i$ and $\left(\Sigma \Sigma^T\right)_{i j }= a_{i j} $ are the drift vector and the diffusion coefficient in \eqref{A:EtKu}.
\end{theorem}

We are actually going to apply this Theorem to a case in which the functions $b_i$ and $a_{ij}$ have an explicit, continuous dependence on the time $t$. This generalization is trivial as it amounts to add one dimension to the state space $\R^d$, introducing the deterministic extra variable $Y(t) = t$. Moreover, for next application of this Theorem, the localization given by the stopping times $\tau^N_h$ will not be necessary, and it will be omitted.

We set
\[   \begin{split}
	X^N (t) \coloneqq 
	\begin{bmatrix}
         \xi^N(t)\\
	\eta^N(t) 
	\end{bmatrix} \coloneqq
	\begin{bmatrix}
	\sqrt{N}\left(m_N(t) - m(t)\right)\\
	\sqrt{N}\left(v_N(t) - v(t)\right)
	\end{bmatrix}
	\end{split}
\]
By using \eqref{eq:micro:dyn2}, \eqref{eq:m(t)v(t)}, the identities $m_N(t) = \frac{\xi^N(t)}{\sqrt{N}}+m(t)$ and $v_N(t) = \frac{\eta^N(t)}{\sqrt{N}}+v(t)$, and recalling the compensated Poisson random measures
defined in \eqref{compPRM}
we have that
\begin{equation}\label{eq:fluctuation:proc}
    \begin{split}
        \xi^N(t) & 
		 =  \int_0^t \left[-(\lambda v(s) + r) \xi^N(s) + \lambda(1-m(s)) \eta^N(s) - \lambda \frac{\xi^N(s) \eta^N(s)}{\sqrt{N}} \right] ds + M^N_1(t)\\
		& \\
		\eta^N(t) & =  \int_0^t \left[ - \lambda v(s) \xi^N(s) + \left( -(\alpha + r) + \lambda (1-m(s)) \right) \eta^N(s) - \lambda \frac{\xi^N(s) \eta^N(s)}{\sqrt{N}} \right] ds + M^N_2(t) 
    \end{split}
\end{equation} 
where 
\begin{equation} 
	\begin{split}
		M^N_1(t)  \coloneqq \xi^N(0) + \sum_{i=1}^N & \left[  \int_{[0,t]\times [0,+\infty)} \frac{\pmb{1}_{[0, (1-\sigma_i(s^-))\lambda v_N(s^-)]}(u)}{\sqrt{N}} \tilde N^\uparrow_i(ds,du) \right. \\ &  - \left.\int_{[0,t]\times [0,+\infty)} \frac{\pmb{1}_{[0,\sigma_i(s^-) r]}(u)}{\sqrt{N}} \tilde N^\downarrow_i(ds,du) \right] 
		\\
		M^N_2(t)  \coloneqq \eta^N(0) + \sum_{i=1}^N  & \left[\int_{[0,t]\times [0,+\infty)} \frac{\pmb{1}_{[0,(1-\sigma_i(s^-)) \lambda v_N(s^-)]}(u)}{\sqrt{N}}\tilde N^\uparrow_i(ds,du) \right. \\ &  - \left. \int_{[0,t]\times [0,+\infty)} \frac{x_i(s^-) \pmb{1}_{[0,\sigma_i(s^-) r]}(u)}{\sqrt{N}}\tilde N^\downarrow_i(ds,du)  \right]
	\end{split}
\end{equation}
are square-integrable martingales having predictable quadratic variation equal to 
\begin{equation}
	\begin{split}
		\langle M^N_1\rangle_t  &= \sum_{i=1}^N \int_{[0,t]\times [0,+\infty)} \left( \frac{\pmb{1}_{[0, (1-\sigma_i(s))\lambda v_N(s)]}(u)}{\sqrt{N}}\right)^2 du ds + \sum_{i=1}^N \int_{[0,t]\times [0,+\infty)}  \left( \frac{\pmb{1}_{[0,\sigma_i(s) r]}(u)}{\sqrt{N}} \right)^2 du ds \\
		& =  \int_0^t \left(\lambda (1-m_N(s)) v_N(s) + r m_N(s) \right) ds\\
		\langle M^N_2\rangle_t & = \sum_{i=1}^N \int_{[0,t]\times [0,+\infty)} \left( \frac{\pmb{1}_{[0, (1-\sigma_i(s))\lambda v_N(s)]}(u)}{\sqrt{N}}\right)^2 du ds + \sum_{i=1}^N \int_{[0,t]\times [0,+\infty)} \left( \frac{x_i(s) \pmb{1}_{[0,\sigma_i(s) r]}(u)}{\sqrt{N}}\right)^2 du ds \\
		& =  \int_0^t \left(\lambda (1-m_N(s)) v_N(s) + r \frac{\sum_{i=1}^N x^2_i(s)}{N} \right) ds.
	\end{split}
\end{equation}
Also, notice that the quadratic covariation process between $M^N_1$ and $M^N_2$ is given by
\begin{equation}
	\begin{split}
	\langle M^N_1, M^N_2\rangle_t & = \sum_{i=1}^N \int_{[0,t]\times [0,+\infty)} \left( \frac{\pmb{1}_{[0, (1-\sigma_i(s))\lambda v_N(s)]}(u)}{\sqrt{N}}\right)^2 du ds + \sum_{i=1}^N \int_{[0,t]\times [0,+\infty)} x_i(s) \left(\frac{\pmb{1}_{[0,\sigma_i(s) r]}(u)}{\sqrt{N}}\right)^2 du ds \\
		& = \int_0^t \left(\lambda (1-m_N(s)) v_N(s) + r v_N(s) \right) ds
	\end{split}
\end{equation}

We can therefore apply Theorem \ref{thm:diffusion:approx} with the following positions:
\[
B^N(t) \coloneqq 
	\begin{bmatrix}   \int_0^t \left[-(\lambda v(s) + r) \xi^N(s) + \lambda(1-m(s)) \eta^N(s) - \lambda \frac{\xi^N(s) \eta^N(s)}{\sqrt{N}} \right] ds \\ \int_0^t \left[ - \lambda v(s) \xi^N(s) + \left( -(\alpha + r) + \lambda (1-m(s)) \right) \eta^N(s) - \lambda \frac{\xi^N(s) \eta^N(s)}{\sqrt{N}} \right] ds 
	\end{bmatrix}
\]
\[
b(x,t) = b(\xi,\eta,t) \coloneqq 
	\begin{bmatrix} -(\lambda v(t) + r) \xi+ \lambda(1-m(t)) \eta \\ - \lambda v(t) \xi + \left( -(\alpha + r) + \lambda (1-m(t)) \right) \eta \end{bmatrix}
\]
\[
A^N(t)  \coloneqq 
	\begin{bmatrix} \int_0^t \left(\lambda (1-m_N(s)) v_N(s) + r m_N(s) \right) ds & \int_0^t \left(\lambda (1-m_N(s)) v_N(s) + r v_N(s) \right) ds \\
	\int_0^t \left(\lambda (1-m_N(s)) v_N(s) + r v_N(s) \right) ds & \int_0^t \left(\lambda (1-m_N(s)) v_N(s) + r \frac{\sum_{i=1}^N x^2_i(s)}{N} \right) ds
	\end{bmatrix}
\]
\[
a(x,t) = a(t) = \begin{bmatrix}
\l \left(1-m(t)\right)v(t) +r m(t) & \l \left(1-m(t)\right)v(t) +r \l v(t) \\
\l\left(1-m(t)\right)v(t) +r v(t) & \l \left(1-m(t)\right) v(t) + r  \E\left[ \bar x^2(t)\right]
\end{bmatrix}
\]
Conditions \eqref{eq:hyp4.4} and \eqref{eq:hyp4.5} are obvious, as $B^N(t)$ and $A^N(t)$ are continuous. Condition \eqref{eq:hyp4.3} is also simple: the jumps of $\xi^N$ and $\eta^N$ are bounded by $\frac{1}{\sqrt{N}}$. We now verify conditions \eqref{eq:hyp4.7} and \eqref{eq:hyp4.6}. We begin by proving  \eqref{eq:hyp4.6} for the first component of $B^N$.
\begin{equation}
    \begin{split}
       \sup_{t\in [0,T]}\Bigg\vert B^N_1(t)-\int_0^t b_1(X^N(s))ds \Bigg\vert 
       & = \sup_{t\in [0,T]}\Bigg\vert \int_0^t \left[ -(\lambda v(s) +r)\xi^N(s) + \lambda (1-m(s)) \eta^N(s) - \lambda \frac{\xi^N(s) \eta^N(s)}{\sqrt{N}} \right. \\
	& \left. + (\lambda v(s) + r) \xi^N(s) -\lambda (1-m(s)) \eta^N(s)
	 \right] ds \Bigg\vert\\
	& \leq T \lambda \sup_{s\in [0,T]}\left| \frac{\xi^N(s) \eta^N(s)}{\sqrt{N}}\right|  
    \end{split}
\end{equation}
Thus, using Corollary \ref{cor:LLN} and Markov inequality, for every $\epsilon> 0$, 
\begin{equation}
    \begin{split}
        & \mathbb{P}\left(\sup_{t \in [0,T]}\Bigg\vert B^N_1(t)-\int_0^t b_1(X^N(s))ds \Bigg\vert \geq \epsilon \right) \leq \mathbb{P}\left(T \lambda \sup_{t\in [0,T]}\left| \frac{\xi^N(t) \eta^N(t)}{\sqrt{N}}\right|   \geq \epsilon
        \right)\\
        &  = \mathbb{P}\left(\sup_{t\in[0,T]}\vert \xi^N(t)\vert\,\sup_{t \in [0,T]} \vert \eta^N(t)\vert\geq \frac{\epsilon\sqrt{N}}{T\lambda}\right)\\
	& \leq \mathbb{P}\left(\sup_{t \in [0,T]} \vert \xi^N(t)\vert \geq \sqrt{\frac{\epsilon}{T\lambda}}N^{1/4}\right)+ \mathbb{P}\left(\sup_{t \in [0,T]} \vert \eta^N(t)\vert \geq \sqrt{\frac{\epsilon}{T\lambda}} N^{1/4}\right)\\
        & = \mathbb{P}\left(\sup_{t\in [0,T]}\vert m_N(t) - m(t)\vert \geq \sqrt{\frac{\epsilon}{T \lambda}}N^{-1/4}\right) + \mathbb{P}\left(\sup_{t\in [0,T]}\vert v_N(t) - v(t)\vert \geq \sqrt{\frac{\epsilon}{T \lambda}}N^{-1/4}\right)\\
        & \leq  2 \frac{\frac{B_T}{\sqrt{N}}}{\sqrt{\frac{\epsilon}{T\lambda}}N^{-1/4}} \overset{N\to+\infty}{\longrightarrow} 0
    \end{split}
\end{equation}
The proof for $B^N_2$ is similar, and it is omitted. We now verify \eqref{eq:hyp4.7} for $i=j=2$, all other cases being similar. 
\begin{equation} \label{an}
    \begin{split}
        & \Bigg\vert A^N_{22} (t) -\int_0^t a_{22}(s) ds \Bigg\vert  =
\Bigg\vert \int_0^t \left( \lambda (1 - m_N(s))v_N(s) + r\frac{\sum_{i=1}^N x^2_i(s)}{N}  - \lambda (1-m(s)) v(s) - r \E\left[ \bar x^2(t)\right]
 \right) ds \Bigg\vert\\
	& \leq \int_0^t \left| \lambda (1-m_N(s)) v_N(s) - \lambda (1-m(s)) v_N(s)\right| ds + \int_0^t \left| \lambda (1-m(s)) v_N(s) - \lambda (1-m(s)) v(s)\right| ds \\
	& + r \int_0^t \left|\frac{\sum_{i=1}^N  x^2_i(s)}{N} -\E\left[ \bar x^2(t)\right]
 \right| ds \\
	& \leq \l \int_0^t\left(  |m_N(s) - m(s) | +  |v_N(s) - v(s)|\right) ds  +  r \int_0^t \left|\frac{\sum_{i=1}^N  x^2_i(s)}{N} -\E\left[ \bar x^2(t)\right] \right| ds
	    \end{split}
\end{equation}
By Markov inequality, we are left to show that 
\[
\lim_{N \ra +\infty}\E\left[ \sup_{t \in [0,T]} \Bigg\vert A^N_{22} (t) -\int_0^t a_{22}(s) ds \Bigg\vert \right]= 0.
\]
Using \eqref{an}, we have
\begin{multline*}
\E\left[ \sup_{t \in [0,T]} \Bigg\vert A^N_{22} (t) -\int_0^t a_{22}(s) ds \Bigg\vert \right] \\ \leq \E\left[\l \int_0^T\left(  |m_N(s) - m(s) | +  |v_N(s) - v(s)|\right) ds  +  r \int_0^T \left|\frac{\sum_{i=1}^N  x^2_i(s)}{N} -\E\left[ \bar x^2(t)\right] \right| \right]
\end{multline*}
which converges to zero as $N \ra +\infty$ by  Corollary \ref{cor:LLN}. 	
\subsection{Proof of Theorem \ref{th:rescal}}

This proof is similar to that of Theorem \ref{thm:fluct:convergence}. However, since some of the parameters are sent to infinity with $N$, we need to use the details of the upper bound for the propagation of chaos (see \eqref{eq:bound:prop:chaos} and \eqref{C_T}). Recall that
\begin{equation}\label{eq:def:fluctuations:rescaled}
    \begin{split}
        &\hat{\xi}^N (t) \coloneqq \sqrt{N}\frac{m_N\left(\frac{t}{\sqrt{\l_N}}\right) - m^\ast}{\l_N^{1/4} } \\
        & \hat{\eta}^N (t) \coloneqq \l_N^{1/4} \sqrt{N}\left(v_N\left(\frac{t}{\sqrt{\l_N}}\right) - v^{\ast} \right).    \end{split}
\end{equation}
We use here the standard scaling invariance of Poisson random measures in the following form.  Let $N(dt,du)$ be a Poisson random measure of intensity $dt \, du$, and let $\L:[0, +\infty) \ra [0,+\infty)$ be right continuous. Then for $a>0$
\[
\int_{[0,at]\times \R^+} \pmb{1}_{[0,\L(s^-)]} N(ds,du) = \int_{[0,t]} \pmb{1}_{[0,a\L(a s^-)]} \hat{N}(ds,du),\
\]
where $\hat{N}$ has the same distribution as $N$. This is checked by defining (interpreting point random measures as random sets)
\[
\hat{N} := \left\{ (t,u): (t/a, au) \in N\right\}.
\]
It is elementary to see that $\hat{N}$ is  a Poisson random measure of intensity $dt \, du$.

We use below this property with $a = \frac{1}{\sqrt{\l_N}}$, and we will omit the superscript \ $\hat{}$ \, on the rescaled Poisson processes.

Using \eqref{eq:micro:dyn2} and \eqref{compPRM} we have
\begin{equation}
	\begin{split}
	\hat{\xi}^N (t) &= \hat{\xi}^N (0) + \int_0^t d \hat{\xi}^N(s) = \hat{\xi}^N(0) + \frac{\sqrt{N}}{\l_N^{1/4}} \frac{1}{N}\sum_{i=1}^N \int_0^t  d \sigma_i\left(\frac{s}{\sqrt{\lambda_N}}\right) \\
			   & = \hat{\xi}^N (0) +  \frac{\sqrt{N}}{\l_N^{1/4}} \frac{1}{N}\sum_{i=1}^N \left\{ \int_{[0,t]\times [0,+\infty)} \pmb{1}_{\left[0,\frac{1}{\sqrt{\lambda_N}}\left(1-\sigma_i\left(\frac{s^-}{\sqrt{\lambda_N}}\right)\right)\lambda_N v_N\left( \frac{s^-}{\sqrt{\lambda_N}}\right)\right]}(u) N^\uparrow_i(ds,du) \right.\\
		& \left. - \int_{[0,t]\times [0,+\infty)} \pmb{1}_{\left[0,\frac{1}{\sqrt{\lambda_N}}\sigma_i\left(\frac{s^-}{\sqrt{\lambda_N}}\right) r\right]}(u) N^\downarrow_i(ds,du) \right\} \\
		& = \hat \xi^N(0) + \frac{\sqrt{N}}{\l_N^{1/4}}\int_0^t \frac{1}{\sqrt{\lambda_N}}\left[\lambda_N\left(1-m_N\left(\frac{s}{\sqrt{\lambda_N}}\right)\right) v_N\left(\frac{s}{\sqrt{\lambda_N}}\right) - r m_N\left(\frac{s}{\sqrt{\lambda_N}}\right) \right] ds \\
		& + \frac{\sqrt{N}}{\l_N^{1/4}}\frac{1}{N}\sum_{i=1}^N \left\{\int_{[0,t]\times [0,+\infty)} \pmb{1}_{\left[0,\frac{1}{\sqrt{\lambda_N}}\left(1-\sigma_i\left(\frac{s^-}{\sqrt{\lambda_N}}\right)\right)\lambda_N v_N\left( \frac{s^-}{\sqrt{\lambda_N}}\right)\right]}(u) \tilde N^\uparrow_i(ds,du) \right.\\
		& \left.  - \int_{[0,t]\times [0,+\infty)} \pmb{1}_{\left[0,\frac{1}{\sqrt{\lambda_N}}\sigma_i\left(\frac{s^-}{\sqrt{\lambda_N}}\right) r\right]}(u) \tilde N^\downarrow_i(ds,du) \right\} \\
		& \eqqcolon  \hat{\xi}^N (0) + \hat B^N_1(t) + \hat M^N_1(t)
\end{split}
\end{equation}
where  
\begin{equation}\label{eq:hatB1}
	\begin{split}
		\hat B^N_1(t) & \coloneqq \frac{\sqrt{N}}{\l_N^{1/4}}\int_0^t \frac{1}{\sqrt{\lambda_N}}\left[\lambda_N\left(1-m_N\left(\frac{s}{\sqrt{\lambda_N}}\right)\right) v_N\left(\frac{s}{\sqrt{\lambda_N}}\right) - r m_N\left(\frac{s}{\sqrt{\lambda_N}}\right) \right] ds \\
		& = \frac{\sqrt{N}}{\l_N^{1/4}}\int_0^t \frac{1}{\sqrt{\lambda_N}}\left[\lambda_N\left(1- m^\ast - \frac{\l_N^{1/4}}{\sqrt{N}}\hat\xi^N(s)\right) \left(v^\ast + \frac{1}{\sqrt{N}\l_N^{1/4}} \hat\eta^N(s)\right) \right.\\
		& \left. - r \left(m^\ast + \frac{\l_N^{1/4}}{\sqrt{N}} \hat\xi^N(s)\right) \right] ds \\
		& \overset{(i)}{=} \int_0^t \frac{1}{\sqrt{\lambda_N}} \left[(-\lambda_N v^\ast - r) \hat\xi^N(s) + \left(\lambda_N (1-m^\ast)\frac{1}{\sqrt{\l_N}}\right)\hat\eta^N(s) -\frac{\lambda_N}{\sqrt{N}\l_N^{1/4}}\hat\xi^N(s) \hat\eta^N(s) \right] ds \\
		& \overset{(ii)}{=} \int_0^t \frac{1}{\sqrt{\lambda_N}}\left[ -\left(\frac{r}{\rho}+ o_N(1)\right) \hat\xi^N(s) + \sqrt{\lambda_N} \left(\rho + o_N(1) \right)\hat\eta^N(s) - \frac{\lambda_N^{3/4}}{\sqrt{N}}\hat\xi^N(s) \hat\eta^N(s) \right] ds, 
	\end{split}
\end{equation}
where in $(i)$ we employed Eq.s \eqref{eq:m(t)v(t)} and in $(ii)$ we employed \eqref{eq:mstar:vstar}, and
\begin{equation}
	\begin{split}
		\hat M^N_1(t) & \coloneqq \hat\xi^N(0) + \frac{\sqrt{N}}{\l_N^{1/4}} \frac{1}{N}\sum_{i=1}^N \left\{\int_{[0,t]\times [0,+\infty)} \pmb{1}_{\left[0,\frac{1}{\sqrt{\lambda_N}}\left(1-\sigma_i\left(\frac{s^-}{\sqrt{\lambda_N}}\right)\right)\lambda_N v_N\left( \frac{s^-}{\sqrt{\lambda_N}}\right)\right]}(u) \tilde N^\uparrow_i(ds,du) \right.\\
		& \left.  - \int_{[0,t]\times [0,+\infty)} \pmb{1}_{\left[0,\frac{1}{\sqrt{\lambda_N}}\sigma_i\left(\frac{s^-}{\sqrt{\lambda_N}}\right) r\right]}(u) \tilde N^\downarrow_i(ds,du) \right\}
	\end{split}
\end{equation}
is a martingale with quadratic variation equal to
\begin{equation}\label{eq:hatM1:var}
	\begin{split}
		\langle \hat M^N_1\rangle_t & = \frac{1}{N\sqrt{\l_N}}\sum_{i=1}^N \int_0^t \frac{1}{\sqrt{\lambda_N}}\left[\lambda_N \left(1-\sigma_i\left(\frac{s}{\sqrt{\lambda_N}}\right)\right) v_N\left(\frac{s}{\sqrt{\lambda_N}}\right) + r \sigma_i\left(\frac{s}{\sqrt{\lambda_N}}\right)\right] ds \\
		& = \frac{1}{\l_N}\int_0^t \left[\lambda_N \left(1-m_N\left(\frac{s}{\sqrt{\lambda_N}}\right)\right) v_N\left(\frac{s}{\sqrt{\lambda_N}}\right) + r m_N\left(\frac{s}{\sqrt{\lambda_N}}\right)\right] ds \\
		& = \int_0^t \left[ \frac{1}{\l_N} \left(\lambda_N(1-m^\ast) v^\ast + r m^\ast \right) + \frac{1}{\sqrt{N}\lambda_N^{3/4}} (-\lambda_N v^\ast  + r )\hat\xi^N(s) \right.\\
		& \left. + \lambda_N (1-m^\ast) \frac{1}{\sqrt{N}\lambda_N^{5/4} } \hat\eta^N(s) -  \frac{\sqrt{\lambda_N} }{N } \hat\xi^N(s) \hat\eta^N(s)\right]ds \\
		& = \int_0^t \left[2 r \left(1-\rho+ o_N(1)\right)\frac{1}{\lambda_N} + \left( -\frac{r}{\rho} + 2 r + o_N(1)\right)  \frac{1}{\sqrt{N}\lambda_N^{3/4}}\hat\xi^N(s) \right.\\
 		& \left. +\left(\rho+ o_N(1) \right) \frac{1}{\sqrt{N} \l_N^{1/4}} \hat\eta^N(s)  - \frac{\sqrt{\lambda_N}}{N } \hat\xi^N(s) \hat\eta^N(s) \right] ds.
	\end{split}
\end{equation}
In \eqref{eq:hatB1} and \eqref{eq:hatM1:var} we denoted by $o_N(1)$ any sequence of real numbers that goes to zero as $N \ra +\infty$. In a similar way we can write
\[
\hat{\eta}^N (t) = \hat{\eta}^N (0) + \hat B^N_2(t) + \hat M^N_2(t),
\]
where
\begin{equation}\label{eq:hatB2}
	\begin{split}
		\hat B^N_2(t) & \coloneqq  \int_0^t \left[ - \lambda_N v^\ast  \hat\xi^N(s) - \frac{\lambda_N^{3/4}}{\sqrt{N}} \hat\xi^N(s) \hat\eta^N(s) \right] ds\\
		& = \int_0^t \left[ - r \left(\frac{1}{\rho} -1 + o_N(1)\right) \hat\xi^N(s) - \frac{\lambda_N^{3/4}}{\sqrt{N}} \hat\xi^N(s) \hat\eta^N(s) \right] ds
	\end{split}
\end{equation}
\begin{equation}
	\begin{split}
		\hat M^N_2(t) & \coloneqq \hat\eta^N(0) + \sqrt{N}\lambda_N^{1/4} \frac{1}{N}\sum_{i=1}^N \left\{\int_{[0,t]\times [0,+\infty)} \pmb{1}_{\left[0,\frac{1}{\sqrt{\lambda_N}}\left(1-\sigma_i\left(\frac{s^-}{\sqrt{\lambda_N}}\right)\right)\lambda_N v_N\left( \frac{s^-}{\sqrt{\lambda_N}}\right)\right]}(u) \tilde N^\uparrow_i(ds,du) \right.\\
		& \left. -  \int_{[0,t]\times [0,+\infty)} x_i\left(\frac{s^-}{\sqrt{\lambda_N}}\right) \pmb{1}_{\left[0,\frac{1}{\sqrt{\lambda_N}}\sigma_i\left(\frac{s^-}{\sqrt{\lambda_N}}\right) r\right]}(u) \tilde N^\downarrow_i(ds,du) \right\}
	\end{split}
\end{equation}
is the remaining martingale after compensation of the PRM's and has quadratic variation equal to 
\begin{equation}\label{eq:hatM2:var}
	\begin{split}
		\langle \hat M^N_2\rangle_t & =  \int_0^t  \left[\lambda_N \left( 1 - m_N\left(\frac{s}{\sqrt{\lambda_N}}\right)\right)v_N\left(\frac{s}{\sqrt{\lambda_N}}\right) + r \frac{1}{N}\sum_{i=1}^N x^2_i\left(\frac{s}{\sqrt{\lambda_N}}\right) \right] ds \\
		& = \int_0^t \left[  \left( \lambda_N (1-m^\ast) v^\ast + r \frac{1}{N}\sum_{i=1}^N x^2_i\left(\frac{s}{\sqrt{\lambda_N}}\right) \right) - \lambda_N v^\ast \frac{\l_N^{1/4}}{\sqrt{N}} \hat\xi^N(s) \right.\\
		& \left. + (1-m^\ast) \frac{\lambda_N^{3/4}  }{\sqrt{N}}\hat\eta^N(s) - \frac{\lambda_N }{N }\hat\xi^N(s) \hat\eta^N(s)\right] ds \\
		& = \int_0^t \left[ \left(r (1-\rho) + r \frac{1}{N}\sum_{i=1}^N x^2_i\left(\frac{s}{\sqrt{\lambda_N}}\right) + o_N(1) \right) - r\left(\frac{1}{\rho}-1 +o_N(1)\right) \frac{\lambda_N^{1/4}}{\sqrt{N}} \hat\xi^N(s) \right.\\
		& \left. + (\rho + o_N(1)) \frac{\lambda_N^{3/4} }{\sqrt{N}}\hat\eta^N(s) - \frac{\lambda_N}{N }\hat\xi^N(s) \hat\eta^N(s)\right] ds
	\end{split}
\end{equation}
We can also compute the covariation of $M^N_1$ and $M^N_2$:
\begin{equation}\label{eq:hat:cov}
	\begin{split}
		\langle \hat M^N_1, \hat M^N_2 \rangle_t & = \sum_{i=1}^N\int_0^t  \frac{1}{N \sqrt{\lambda_N}} \left[  \left( 1 - \sigma_i\left(\frac{s}{\sqrt{\lambda_N}}\right)\right) \lambda_N v_N\left(\frac{s}{\sqrt{\lambda_N}}\right) + r x_i\left(\frac{s}{\sqrt{\lambda_N}}\right) \right]  ds \\
		& = \frac{1}{ \sqrt{\lambda_N}}  \int_0^t \left[\left(1-m^\ast - \frac{\l_N^{1/4}}{\sqrt{N}}\hat\xi^N(s)\right)\lambda_N \left(v^\ast + \frac{1}{\sqrt{N}\lambda_N^{1/4}} \hat\eta^N(s)\right) \right.\\
		& \left. + r\left(v^\ast + \frac{1}{\sqrt{N}\lambda_N^{1/4}} \hat\eta^N(s)\right) \right] ds  \\
		& =  \frac{1}{ \sqrt{\lambda_N}}  \int_0^t \left[ r(1-\rho + o_N(1)) + \frac{r^2}{\lambda_N}\left(\frac{1}{\rho}-1 + o_N(1)\right) - r \left(\frac{1}{\rho}-1 + o_N(1)\right) \frac{\l_N^{1/4}}{\sqrt{N}} \hat\xi^N(s) \right.\\
		& \left. + (\lambda_N (\rho+o_N(1))+ r) \frac{1}{\sqrt{N}\lambda_N^{1/4}} \hat\eta^N(s) - \frac{1}{N } \hat\xi^N(s) \hat\eta^N(s)\right] ds.
	\end{split}
\end{equation}
We can therefore apply Theorem \ref{thm:diffusion:approx} with
\[
b(x,t) = b(\xi,\eta) \coloneqq 
	\begin{bmatrix} \rho\eta \\ - \frac{r}{\rho}(1-\rho) \xi \end{bmatrix}
\]
\[
a(x,t) \coloneqq 
	\begin{bmatrix} 0 & 0 \\ 0  & r(1-\rho) \end{bmatrix}.
\]
Using the localization by the stopping time $\tau_h^N$ defined in Theorem \ref{thm:diffusion:approx}, all conditions required by Theorem \ref{thm:diffusion:approx} are readily checked except for the convergence of $\langle \hat M^N_2\rangle_t$, where the key point is to show that for $\e>0$
\be{v2cond}
\lim_{N \ra +\infty} \P\left(\sup_{0 \leq t \leq \frac{T}{\sqrt{\l_N}}} \frac{1}{N}\sum_{i=1}^N x^2_i\left(t\right) > \e \right) = 0,
\ee
which in turn follows from
\be{condv2}
\lim_{N \ra +\infty} \E\left[\sup_{0 \leq t \leq \frac{T}{\sqrt{\l_N}}} \frac{1}{N}\sum_{i=1}^N x^2_i\left(t\right)\right] = 0.
\ee
By using \eqref{C_Test2} and observing that 
\[
\E[\bar{x}^2(t)] \leq \E[\bar{x}(t)] = v^* = \frac{r}{r+\a_N} \left(1-\frac{r+\a_N}{\l_N} \right) \ra 0
\]
as $N \ra +\infty$, to prove \eqref{v2cond} it is enough to show that if $C_T$ is the time-dependent constant in  \eqref{CT3}, we have
\be{CTlim}
\lim_{N \ra +\infty} \frac{C_{T/\sqrt{\l_N}}}{\sqrt{N}} = 0.
\ee
Assuming
\[
\lim_{N \ra +\infty} \l_N = +\infty,
\]
it is easily checked that the dominant term in $\frac{C_{T/\sqrt{\l_N}}}{\sqrt{N}}$ is of order
\[
\frac{\sqrt{\l_N}}{\sqrt{N}} \exp\left[ KT \sqrt{\l_N} \right] = \exp\left[KT \sqrt{\l_N}  + \frac12 \log \l_N - \frac12 \log N \right],
\]
for a suitable constant $K$, which goes to zero provided 
\[
\lim_{N \ra +\infty} \frac{\sqrt{\l_N}}{\log N} = 0.
\]
This completes the proof.

\medskip
\noindent
{\bf Acknowledgments.} We are deeply grateful to Francesca Collet, Marco Formentin and Markus Fischer for several key suggestions, that had a strong impact in the development of this research. EM acknowledges financial support from Progetto Dottorati—Fondazione Cassa di Risparmio di Padova e Rovigo. 

\bibliographystyle{abbrv}

\begin{thebibliography}{99}

\bibitem{DPMarsiglia} Dai Pra, Paolo. \lq\lq Stochastic mean-field dynamics and applications to life sciences.\rq\rq\, International workshop on Stochastic Dynamics out of Equilibrium. Cham: Springer International Publishing, 2017.
\bibitem{gardiner} Gardiner, Crispin W. Handbook of stochastic methods - for physics, chemistry and the natural sciences, Second Edition, Springer Series in Synergetics, 1986. 
\bibitem{G92} Graham, Carl. \lq\lq Nonlinear diffusion with jumps.\rq\rq\, Annales de l'IHP Probabilités et statistiques. Vol. 28. No. 3. 1992.
\bibitem{L13} Liggett, Thomas M. Stochastic interacting systems: contact, voter and exclusion processes. Vol. 324. Springer science \& Business Media, 2013.

\bibitem{EtKu2005}
Ethier, S. N. and Kurtz, T. G. 
\newblock Markov Processes:  Characterization and Convergence.  
\newblock 2005.

\bibitem{Ikeda}
Ikeda, N. and Watanabe, S. 
\newblock Stochastic Differential Equations and Diffusion Processes. 
\newblock 2014.

\bibitem{MKN}
McKane, A. J. and Newman, T. J. \lq\lq Predator-Prey Cycles from Resonant Amplification of Demographic Stochasticity.\rq\rq\, Phys. Rev. Lett. Vol. 94. No. 21. Pages 218102. 2005.

\bibitem{ACF}
Albi, G., Chignola, R. and Ferrarese, F. \lq\lq Efficient ensemble stochastic algorithms for agent-based models with spatial predator–prey dynamics.\rq\rq\, Mathematics and Computers in Simulation,
Vol. 199. Pages 317-340. 2022.



\end{thebibliography}

\end{document}